\documentclass[12pt]{amsart}
\usepackage{amsfonts}

\usepackage{float}
\usepackage{graphicx}
\usepackage{amscd}
\usepackage{amsmath}

\theoremstyle{plain}
\newtheorem{X}{X}[section]
\newtheorem{theorem}[X]{Theorem}
\newtheorem*{theorem*}{Theorem}
\newtheorem{proposition}[X]{Proposition}
\newtheorem*{proposition*}{Proposition}
\newtheorem{lemma}[X]{Lemma}
\newtheorem*{lemma*}{Lemma}
\newtheorem{corollary}[X]{Corollary}

\theoremstyle{definition}

\newtheorem{remark}[X]{Remark}
\newtheorem*{remark*}{Remark}

\newtheorem*{example*}{Example}

\numberwithin{equation}{section}

\usepackage{verbatim}
\usepackage[size=2em,heads=vee]{diagrams}
\diagramstyle[labelstyle=\scriptstyle]

\makeindex

\newenvironment{mylist}
{\begin{list}
{--}
{\setlength{\leftmargin}{.5in}\setlength{\rightmargin}{.5in}}}
{\end{list}}

\numberwithin{equation}{section}

\def\func#1{\mathop{\rm #1}\nolimits}%
\def\1{{1\mkern-7mu1}}  

\newcommand\End{\operatorname{End}}

\newcommand\Gal{\operatorname{Gal}}
\newcommand\GL{\operatorname{GL}}

\newcommand\Hom{\operatorname{Hom}}

\newcommand\inv{\operatorname{inv}}

\newcommand\Ker{\operatorname{Ker}}

\newcommand\Nm{\operatorname{Nm}}

\newcommand\ord{\operatorname{ord}}

\newcommand\Res{\operatorname{Res}}

\newcommand\Tr{\operatorname{Tr}}

\def\al{{ \text{al} }}

\def\ord{{\text{ord}}}

\begin{document}
\title[Tate Conjecture]{The Tate Conjecture for Certain Abelian Varieties over Finite Fields}
\author{J.S. MILNE}
\address{Mathematics Department, University of Michigan, Ann Arbor, MI 48109.}
\email{jmilne@umich.edu}
\dedicatory{August 1, 1999; submitted version.}
\thanks{This article includes research supported in part by the National Science
Foundation}
\maketitle

In an earlier work, I showed that if the Hodge conjecture holds for all\emph{%
\ } complex abelian varieties of CM-type, then the Tate conjecture holds for
all\emph{\/} abelian varieties over finite fields (Milne 1999b). In this
article, I extract from the proof a statement (Theorem 1.1) that sometimes
allows one to deduce the Tate conjecture for the powers of single abelian
variety $A$ over a finite field from knowing that some Hodge classes on
their lifts to characteristic zero are algebraic.

Tate's theorem (Tate 1966) implies that the Tate conjecture holds for any
abelian variety over a finite field whose $\mathbb{Q}_{\ell }$-algebra of
Tate classes is generated by the classes of degree $1$. Using Theorem 1.1
and a result of Schoen (1988), I construct families of abelian varieties
over finite fields for which this condition fails, but which nevertheless
satisfy Tate's conjecture (Theorem 1.7).

The main results are stated in Section 1 and proved in Section 2. Appendix A
summarizes the theories of abelian varieties of CM-type over $\mathbb{C}$
and of abelian varieties over finite fields, and how the reduction map
relates the two. Appendix B sharpens a result of Clozel on the relation
between numerical and homological equivalence for abelian varieties over
finite fields.

Notations not introduced in \S1 are listed at the start of Appendix A.

\section{Statements}

Let $X$ be a smooth complete variety over an algebraic closure $\mathbb{F}$
of the field $\mathbb{F}_{p}$ of $p$ elements. The choice of a model $X_{1}$
of $X$ over a subfield $\mathbb{F}_{p^{n}}$ of $\mathbb{F}$ determines an
action of $\Gal(\mathbb{F}/\mathbb{F}_{p^{n}})$ on the \'{e}tale cohomology
group $H^{2r}(X,\mathbb{Q}_{\ell }(r))$, and we define 
\begin{equation*}
\mathcal{T}_{\ell }^{r}(X)=\bigcup_{X_{1}/\mathbb{F}_{p^{n}}}H^{2r}(X,%
\mathbb{Q}_{\ell }(r))^{\Gal(\mathbb{F}/\mathbb{F}_{p^{n}})}
\end{equation*}
(union over all models). The elements of $\mathcal{T}_{\ell }^{r}(X)$ are
called the $\ell $\emph{-adic Tate classes } \emph{of degree} $r$ on $X$. We
shall say that \emph{the Tate conjecture holds for} $X$ if the $\mathbb{Q}%
_{\ell }$-vector space $\mathcal{T}_{\ell }^{r}(X)$ is spanned by the
classes of algebraic cycles for all $r$ and all $\ell \neq p$.

A Tate class is said to be \emph{exotic\/} if it is not in the $\mathbb{Q}%
_{\ell }$-algebra generated by the Tate classes of degree $1$. For an
abelian variety over $\mathbb{F}{}$, Tate (1966) showed that all Tate
classes of degree $1$ are divisor classes, and so the nonexotic Tate classes
are algebraic.

For an abelian variety $A$ over $\mathbb{F}$, Weil showed that the Frobenius
endomorphism acts semisimply on $H^{1}(A,\mathbb{Q}_{\ell })$, and hence on
all the \'{e}tale cohomology groups of $A$. It follows that, if the Tate
conjecture holds for $A$, then the statements of Tate 1994, Theorem 2.9,
hold for every model $A_{1}$ of $A$ over a finite field. In particular, the
Tate conjecture holds for $A_{1}/\mathbb{\mathbb{F}{}}_{q}$, and for every $%
r $, the order of the pole of the zeta function $Z(A_{1},t)$ of $A_{1}$ at $%
t=q^{-r}$ is equal to the rank of the group of numerical equivalence classes
of algebraic cycles of codimension $r$ on $A_{1}$.

Let $A$ be an abelian variety with many endomorphisms over an algebraically
closed field $k$. Then (see A.3) there is a group of multiplicative type $%
L(A)$ over $\mathbb{Q}$ whose fixed tensors in any Weil cohomology of a power%
$A^{s}$ of $A$ are exactly the \emph{Lefschetz classes}, i.e., those in the
algebra generated by divisor classes. We call $L(A)$ the \emph{Lefschetz
group\/} of $A$.

Now take $k=\mathbb{Q}^{\al}$, and let $w_{0}$ be a prime of $\mathbb{Q}^{\al%
}$ dividing $p$. It follows from the theory of N\'{e}ron models, that $A$
has good reduction at $w_{0}$ (Serre and Tate 1968, Theorem 6), and so
defines an abelian variety $A_{0}$ over the residue field $\mathbb{F}$ at $w$%
. There is a canonical inclusion $L(A_{0})\hookrightarrow L(A)$.

Let $H^{r}(A,\mathbb{Q})$ denote the usual cohomology group of the complex
manifold $A(\mathbb{C})$, and let $H^{r}(A,\mathbb{Q}(m))=H^{r}(A,(2\pi
i)^{m}\mathbb{Q})$ --- it is a rational Hodge structure of weight $r-2m$.
The action of $L(A)$ on $H^{2r}(A^{s},\mathbb{Q}(r))$ defines a
decomposition 
\begin{equation*}
H^{2r}(A^{s},\mathbb{Q}(r))\otimes \mathbb{Q}^{\al}=\bigoplus_{\chi \in
X^{\ast }(L(A))}H^{2r}(A^{s},\mathbb{Q}(r))_{\chi }.
\end{equation*}
We say that $\chi $ is \emph{algebraic\/} if $H^{2r}(A^{s},\mathbb{Q}%
(r))_{\chi }$ contains a nonzero algebraic class for some $r$ and $s$. The
set of algebraic characters of $L(A)$ is stable under the action of $\Gal(%
\mathbb{Q}^{\al}/\mathbb{Q})$, and if $\chi $ is algebraic then $%
H^{2r}(A^{s},\mathbb{Q}(r))_{\chi }$ consists entirely of algebraic classes.%
\footnote{%
The algebraic characters are precisely those that are trivial on the
subgroup $M(A)$ of $L(A)$ --- see A.3.} By composition, an algebraic
character of $L(A)$ defines a character of $L(A_{0})$.

Let $P(A_{0})$ be the smallest algebraic subgroup of $L(A_{0})$ containing a
representative of the Frobenius germ (see A.3).

\begin{theorem}
If 
\begin{equation*}
P(A_{0})=\bigcap \Ker(\chi \colon L(A_{0})\rightarrow \mathbb{G}_{m})
\end{equation*}
(intersection over the algebraic characters of $L(A)$), then the Tate
conjecture holds for all powers of $A_{0}$.
\end{theorem}

An element of $H^{2r}(A,\mathbb{Q}(r))\cap H^{0,0}$ is called a \emph{Hodge
class of degree} $r$ on $A$. We say that \emph{the Hodge conjecture holds for%
} $A$ if the $\mathbb{Q}$-vector space of Hodge classes on $A$ of degree $r$
is spanned by the classes of algebraic cycles for all $r$. The \emph{%
Mumford-Tate group} $MT(A)$ of $A$ is the largest subgroup of $L(A)$ whose
elements fix the Hodge classes on all powers of $A$.

\begin{corollary}
If the Hodge conjecture holds for all powers of $A$ and 
\begin{equation*}
P(A_0)=L(A_0)\cap MT(A)\quad\text{(intersection inside }L(A)) ,
\end{equation*}
then the Tate conjecture holds all powers of $A_0$.
\end{corollary}

A Hodge class is said\footnote{%
Following Tate 1994, p.\ 82.} to be \emph{exotic\/} if it is not in the $%
\mathbb{Q}$-algebra generated by Hodge classes of degree $1$. Lefschetz
showed that all Hodge classes of degree $1$ are divisor classes, and so the
nonexotic Hodge classes are exactly the Lefschetz classes (in particular,
they are algebraic).

Let $E$ be a CM-field of degree $2n$, $n>2$, over $\mathbb{Q}$ containing a
quadratic imaginary field $Q$. Choose an embedding $\rho _{0}\colon
Q\rightarrow \mathbb{Q}^{\al}$, and let $\{\sigma _{0},\ldots ,\sigma
_{n-1}\}$ be the set of extensions of $\rho _{0}$ to $E$. Then $\Phi _{0}=_{%
\text{df}}\{\sigma _{0},\iota \sigma _{1},\ldots ,\iota \sigma _{n-1}\}$ is
a CM-type on $E$ and $\Phi =_{\text{df}}\{\rho _{0}\}$ is a CM-type on $Q$.
Let $(A,i)$ and $(B,j)$ be abelian varieties over $\mathbb{Q}^{\al}$ of
CM-types $(E,\Phi _{0})$ and $(Q,\Phi )$ respectively. We let $Q$ act
diagonally on $A\times B^{n-2}$.

\begin{lemma}
The exotic Hodge classes on $A\times B^{n-2}$ are exactly the nonzero
elements of the subspace 
\begin{equation*}
W(A,B)=_{\text{df}}(\bigwedge_{Q}^{2n-2}H^{1}(A\times B^{n-2},\mathbb{Q}%
))(n-1)
\end{equation*}
of $H^{2n-2}(A\times B^{n-2},\mathbb{Q}(n-1))$.
\end{lemma}

As $A\times B^{n-2}$ has dimension $2n-2$, $H^{1}(A\times B^{n-2},\mathbb{Q}%
) $ has dimension $4n-4$ over $\mathbb{Q}{}$, and so $%
\bigwedge_{Q}^{2n-2}H^{1}(A\times B^{n-2},\mathbb{Q})$ has dimension $1$
over $Q$. The action of an endomorphism of an abelian variety on its
cohomology groups preserves algebraic classes, and so, if $W(A,B)$ contains
one nonzero algebraic class, then it is spanned as a $\mathbb{Q}$-space by
algebraic classes.

\begin{theorem}
If some exotic Hodge class on $A\times B^{n-2}$ is algebraic, then the Hodge
conjecture holds for all abelian varieties of the form $A^{s}\times B^{t}$, $%
s,t\in \mathbb{N}{}$.
\end{theorem}

The abelian varieties $A$ and $B$ over $\mathbb{Q}{}^{\text{al}}$ reduce
modulo $w_{0}$ to abelian varieties $A_{0}$ and $B_{0}$ over $\mathbb{F}{}$.
Let $K$ be the Galois closure of $\sigma _{0}E$ in $\mathbb{Q}{}^{\text{al}}$%
, and let $D(w_{0})$ be the decomposition group of $w_{0}$ in $\Gal(K/%
\mathbb{Q}{})$.

\begin{theorem}
Assume $p$ splits in $Q$ and that $\Gal(K/\sigma _{0}E)\cdot
D(w_{0})=D(w_{0})\cdot \Gal(K/\sigma _{0}E)$.

\begin{enumerate}
\item  For all $\ell \neq p$, the exotic $\ell $-adic Tate classes on $%
A_{0}\times B_{0}^{n-2}$ are exactly the nonzero elements of the subspace 
\begin{equation*}
W(A_{0},B_{0})=_{\text{df}}(\bigwedge_{Q\otimes _{\mathbb{Q}}\mathbb{Q}%
_{\ell }}^{2n-2}H^{1}(A_{0}\times B_{0}^{n-2},\mathbb{Q}_{\ell }))(n-1)
\end{equation*}
of $H^{2n-2}(A_{0}\times B_{0}^{n-2},\mathbb{Q}_{\ell }(n-1))$.

\item  If some exotic Hodge class on $A\times B^{n-2}$ is algebraic, then
the Tate conjecture holds for all abelian varieties over $\mathbb{F}$ of the
form $A_{0}^{s}\times B_{0}^{t}$, $s,t\in \mathbb{N}$.
\end{enumerate}
\end{theorem}

\begin{remark}
(a) Note that, under the hypotheses of the theorem, the $\mathbb{Q}{}$%
-algebra of Hodge classes on $A\times B^{n-2}$ is larger than the tensor
product of the similar algebras for $A$ and $B^{n-2}$, and the $\mathbb{Q}%
{}_{\ell }$-algebra of Tate classes on $A_{0}\times B_{0}^{n-2}$ is larger
than the tensor product of the similar algebras for $A_{0}$ and $B_{0}^{n-2}$%
. Moreover, the groups $L(A\times B)$ and $MT(A\times B)$ (resp. $%
L(A_{0}\times B_{0})$ and $P(A_{0}\times B_{0})$) are not distinguished by
their fixed tensors in the cohomology of $A\times B$ (resp. $A_{0}\times
B_{0}$).

(b) The condition $\Gal(K/E)\cdot D(w_{0})=D(w_{0})\cdot \Gal(K/E)$ holds,
for example, if $E$ is Galois over $\mathbb{\mathbb{Q}{}}$. Without it, the
analysis becomes very complicated, and the theorem fails.
\end{remark}

\subsection{Examples}

Let $C$ be an abelian variety over $\mathbb{C}$, and let $i\colon
Q\rightarrow \End^{0}(C)$ be a homomorphism of $\mathbb{Q}$-algebras, where,
as above, $Q$ is a quadratic imaginary extension of $\mathbb{Q}{}$. The pair 
$(C,i)$ is said to be of \emph{Weil type\/} if the tangent space to $C$ at $%
0 $ is a free $Q\otimes _{\mathbb{Q}}\mathbb{C}$-module.

When $(C,i)$ is of Weil type, its dimension is even, say, $\dim C=2n$, and
the subspace $(\bigwedge_{Q}^{2n}H^{1}(C,\mathbb{Q}))(n)$ of $H^{2n}(C,%
\mathbb{Q}(n))$ consists of Hodge classes (Weil 1977) --- they are called
the \emph{Weil classes\/} on $C$.

Let $\lambda $ be a polarization of $C$ whose Rosati involution induces
complex conjugation on $Q$, and let $E^{\lambda }$ be the Riemann form
defined by $\lambda $. There exists a skew-Hermitian form $\phi \colon
H_{1}(A,\mathbb{Q})\times H_{1}(A,\mathbb{Q})\rightarrow Q$ such that $\Tr%
_{Q/\mathbb{Q}}\circ \phi =E$. The discriminant of $\phi $ is an element of $%
\mathbb{Q}^{\times }/\Nm(Q^{\times })$ which is independent of the choice of
the polarization, and so can be denoted by $\det (C,i)$. The quotient $%
\mathbb{Q}^{\times }/\Nm(Q^{\times })$ is an infinite group killed by $2$,
and for any $a\in \mathbb{Q}^{\times }/\Nm(Q^{\times })$ with $(-1)^{n}a>0$,
there exists an $n^{2}$-dimensional family of abelian varieties of Weil type
with determinant $a$ (Weil 1977, van Geemen 1994).

The following statement is proved in Schoen 1988.

\begin{quote}
If $C$ has dimension $4$, the field $Q$ is generated by a root of $1$ (so, $%
Q $ is $\mathbb{Q}[\sqrt{-3}]$ or $\mathbb{Q}[\sqrt{-1}]$), and $\det
(C,i)=1 $, then the Weil classes on $C$ are algebraic.
\end{quote}

\noindent See also van Geemen 1994, 4.15, 7.1, where it is noted that the
determinant condition was omitted in Schoen 1988, and van Geemen 1996, where
a different proof is given in the case $Q=\mathbb{Q}[\sqrt{-1}]$.

\begin{theorem}
Let $A,B,E,Q$ be as in Theorem 1.5, and assume 
\begin{mylist}
\item  \ $Q$ is generated by a root of $1$;
\item  \ when regarded as an abelian variety of Weil type through the diagonal action
 of $Q$, $A\times B$ has determinant $1$  (modulo squares);
\item  \ $E$ has degree $6$ over $\mathbb{Q}$;
\item \ $p$ splits in $Q$.
\end{mylist}
Then the Hodge conjecture hold for all abelian varieties of the form $%
A^{s}\times B^{t}$, $s,t\in \mathbb{N}\times \mathbb{N}$, and the Tate
conjecture holds for all abelian varieties of the form $A_{0}^{s}\times
B_{0}^{t}$, $s,t\in \mathbb{N}\times \mathbb{N}$.
\end{theorem}

\begin{proof}
The four-dimensional abelian variety $A\times B$ satisfies the hypotheses of
Schoen's theorem, and so Theorem 1.5 applies.
\end{proof}

\begin{remark}
\begin{enumerate}
\item  In the situation of the Theorem 1.7, the space $W(A,B)$ consists of
exotic Hodge classes, and the space $W(A_{0},B_{0})$ consists of exotic Tate
classes.

\item  Similar results hold for $\ell =p$ when the \'{e}tale cohomology
groups are replaced by the crystalline cohomology groups.
\end{enumerate}
\end{remark}

\section{Proofs}

Notations concerning groups of multiplicative type are reviewed at the start
of Appendix A.

\subsubsection{Proof of 1.1}

After the theorem in A.3, in order to prove Theorem 1.1, it suffices to show
that its hypotheses imply that $M(A_{0})=P(A_{0})$.

As numerical equivalence agrees with homological equivalence in
characteristic zero (see B.1), we may regard $M(A)$ as the subgroup of $L(A)$
fixing the algebraic classes in $H^{2r}(A^{s},\mathbb{Q}(r))$ for all $r,s$,
i.e., as the intersection of the kernels of the algebraic characters on $%
L(A) $. Hence 
\begin{equation*}
L(A_{0})\cap M(A)=\bigcap_{\chi \text{ algebraic}}\Ker(\chi \colon
L(A_{0})\rightarrow \mathbb{G}_{m}).
\end{equation*}
Thus, the hypotheses of Theorem 1.1 imply that $L(A_{0})\cap M(A)=P(A_{0})$.
Since 
\begin{equation*}
L(A_{0})\cap M(A)\supset M(A_{0})\supset P(A_{0}),
\end{equation*}
this implies that $M(A_{0})=P(A_{0})$.

\subsubsection{Proof of 1.2.}

As we noted in the proof of 1.1, 
\begin{equation*}
M(A)=\bigcap_{\chi \text{ algebraic}}\Ker(\chi \colon L(A)\rightarrow 
\mathbb{G}_{m}).
\end{equation*}
If the Hodge conjecture holds for the powers of $A$, then $MT(A)=M(A)$ (see
A.3). If, in addition, $P(A_{0})=L(A_{0})\cap MT(A)$, then 
\begin{equation*}
P(A_{0})=\bigcap_{\chi }\Ker(\chi \colon L(A_{0})\rightarrow \mathbb{G}_{m})
\end{equation*}
(intersection over the algebraic characters of $L(A)$), and so (1.2) follows
from (1.1).

\subsubsection{Proofs of 1.3 and 1.4}

Let $E$, $Q$, $\rho _{0}$, $\{\sigma _{0},\ldots ,\sigma _{n-1}\}$, $\Phi $
and $\Phi _{0}$ be as in the paragraph preceding the statement of Lemma 1.3.
Let $K$ be a CM-subfield of $\mathbb{Q}^{\text{al}}$, finite and Galois over 
$\mathbb{Q}{}$, containing the Galois closure of $\sigma _{0}E$ in $\mathbb{Q%
}{}^{\text{al}}$, and let $S^{K}$ be its Serre group (see A.4). For each $i$%
, $0\leq i\leq n-1$, let 
\begin{equation*}
\Sigma _{i}=\{\tau \in \Gal(K/\mathbb{Q}{})\mid \tau \circ \sigma
_{0}=\sigma _{i}\}.
\end{equation*}
Then $\Sigma _{0}$ is the subgroup $\Gal(K/\sigma _{0}E)$ of $\Gal(K/\mathbb{%
Q}{})$ and $\Sigma _{0},\ldots ,\Sigma _{n-1},\iota \Sigma _{0},\ldots
,\iota \Sigma _{n-1}$ are its left cosets. Let $\psi _{i}$ be the
characteristic function of $\Sigma _{i}\cup \bigcup_{j\neq i}\iota \Sigma
_{j}$, and let $\psi $ be the characteristic function of $\bigcup \Sigma
_{i}=\{\tau \mid \tau \circ \rho _{0}=\rho _{0}\}$. Note that $\Sigma _{K}$
acts on the set $\{\Sigma _{0},\ldots ,\iota \Sigma _{n-1}\}$, and that if $%
\tau \Sigma _{i}=\Sigma _{i^{\prime }}$, then $\tau \psi _{i}=\psi
_{i^{\prime }}$. The linear relations among $\psi _{0},\ldots ,\psi
_{n-1},\psi ,\iota \psi $ regarded as elements of $X^{\ast }(S^{K})$ are
exactly the multiples of 
\begin{equation*}
\psi _{0}+\cdots +\psi _{n-1}+(n-2)\psi =(n-1)(\psi +\iota \psi )\qquad 
\text{(*)}.
\end{equation*}

Let $(A,i)$ be an abelian variety of CM-type $(E,\Phi _{0})$, and identify $%
X^{\ast }(L(A))$ with a quotient of $\mathbb{Z}{}^{\Sigma _{E}}$ (see A.5).
The map 
\begin{equation*}
X^{\ast }(\rho _{\Phi _{0}})\colon X^{\ast }(L(A))\rightarrow X^{\ast
}(S^{K})
\end{equation*}
(ibid.) sends $[\sigma _{0}]$ to $\psi _{0}$ and hence, by equivariance and
linearity, it sends $[\sigma _{i}]$ to $\psi _{i}$ and $[\sigma _{0}+\iota
\sigma _{0}]$ to $\psi _{0}+\iota \psi _{0}=\psi +\iota \psi $. Because $%
[\sigma _{0}],\ldots ,[\sigma _{n-1}],[\sigma _{0}+\iota \sigma _{0}]$ form
a basis for $X^{\ast }(L(A))$ and $\psi _{0},\ldots ,\psi _{n-1},\psi +\iota
\psi $ are linearly independent in $X^{\ast }(S^{K})$, we see that $X^{\ast
}(\rho _{\Phi _{0}})\colon X^{\ast }(L(A))\rightarrow X^{\ast }(S^{K})$ is
injective. Therefore, $\rho _{\Phi _{0}}\colon S^{K}\rightarrow L(A)$ is
surjective, and $MT(A)=L(A)$ (ibid.). Hence all Hodge classes on all powers
of $A$ are Lefschetz (A.3, Theorem). In particular, the Hodge conjecture
holds for $A$ and its powers.

Let $(B,i)$ be an elliptic curve of CM-type $(Q,\Phi )$. In this case, $%
L(B)=(\mathbb{\mathbb{G}{}}_{m})_{Q/\mathbb{Q}{}}$ and $X^{\ast }(L(B))=%
\mathbb{Z}{}^{\Sigma _{Q}}$. The map $X^{\ast }(\rho _{\Phi })$ sends $\rho
_{0}$ to $\psi $ and $\iota \rho _{0}$ to $\iota \psi $. As $\psi $ and $%
\iota \psi $ are linearly independent in $X(S^{K})$, this shows that $%
MT(B)=L(B)$, and so all Hodge classes on all powers of $B$ are Lefschetz.

The abelian varieties $A$ and $B$ are simple and nonisogenous (because they
have different dimensions). Their product $A\times B$ is of CM-type $%
(E\times Q,\Phi ^{\prime })$ where $\Phi ^{\prime }=\Phi _{0}\sqcup \Phi $.
The group $X^{\ast }(L(A\times B))$, regarded as a quotient of $\mathbb{Z}%
{}^{\Sigma _{E}\sqcup \Sigma _{Q}},$ has basis $\{[\sigma _{0}],\ldots
,[\sigma _{n-1}],[\rho _{0}],[\rho _{0}+\iota \rho _{0}]\}$, and $X^{\ast
}(\rho _{\Phi ^{\prime }})$ sends 
\begin{equation*}
\lbrack \sigma _{i}]\mapsto \psi _{i},\quad \lbrack \rho _{0}]\mapsto \psi
,\quad \lbrack \rho _{0}+\iota \rho _{0}]\mapsto \psi +\iota \psi .
\end{equation*}
As (*) is the only relation among $\psi _{0},\ldots ,\psi _{n-1},\psi ,\iota
\psi $, the kernel of $X^{\ast }(L(A\times B))\rightarrow X^{\ast }(S^{K})$
is free of rank $1$ with generator 
\begin{equation*}
\chi =[\sigma _{0}+\cdots +\sigma _{n-1}+(n-2)\rho _{0}-(n-1)(\rho
_{0}+\iota \rho _{0})].
\end{equation*}
Hence, there is an exact sequence 
\begin{equation*}
0\rightarrow MT(A\times B)\rightarrow L(A\times B)\overset{\chi }{%
\rightarrow }T\rightarrow 0
\end{equation*}
where $T$ is the $1$-dimensional torus over $\mathbb{Q}$ whose character
group $\langle \chi \rangle $ is isomorphic to $\mathbb{Z}$ with $\Gal(K/%
\mathbb{Q})$ acting (nontrivially) through $\Gal
(Q/\mathbb{Q})$.

The exotic Hodge classes on $A\times B$ and its powers are those that lie in
a rational subspace on which $L(A\times B)$ acts through the characters $%
m\chi $, $m\neq 0$.

We now prove 1.3. The Lefschetz group of $A\times B^{n-2}$ equals that of $%
A\times B$. It acts on 
\begin{equation*}
W(A,B)=\bigwedge_{Q}^{n}H^{1}(A)\otimes
\bigwedge_{Q}^{n-2}((n-2)H^{1}(B))\otimes \mathbb{Q}{}(n-1)
\end{equation*}
through the characters $\chi $ and $\iota \chi =-\chi $. Because $\chi $ is
trivial on $MT(A\times B)$, this space consists of Hodge classes, and
because $\chi $ is not trivial on $L(A\times B)$, the Hodge classes are
exotic. The group $L(A\times B)$ acts on no other subspace of a space $%
H^{2r}(A\times B^{n-2},\mathbb{Q}(r))\otimes \mathbb{Q}^{\al}$ through the
characters $\pm \chi $, and so the elements of $W(A,B)$ are the only exotic
Hodge classes on $A\times B^{n-2}$.

We now prove 1.4. If some exotic Hodge class in $A\times B^{n-2}$ is
algebraic, then $\chi $ is trivial on $M(A\times B)$. Hence $M(A\times
B)=MT(A\times B)$. But $M(A^{s}\times B^{t})=M(A\times B)$ and $%
MT(A^{s}\times B^{t})=MT(A\times B)$ for any $s,t\geq 1$ (see A.5), and so
the Hodge conjecture holds for $A^{s}\times B^{t}$ (see A.3).

\subsubsection{Proof of 1.5.}

We shall compute the terms in the diagram 
\begin{equation*}
\begin{CD} S^K@>>>L(A\times B)\\ @AAA@AAA\\ P^K@>>>L(A_0\times B_0) \end{CD}
\end{equation*}
or, equivalently, in the corresponding diagram of character groups. In fact,
we shall prove that there is an exact commutative diagram 
\begin{equation*}
\begin{CD} 0@>>>\langle\chi\rangle@>>>X^*(L(A\times B))@>>>X^*(S^K)\\
@.@VV{\cong}V@VVV@VVV\\ 0@>>>\langle\chi_{0}\rangle@>>>X^*(L(A_0\times
B_0))@>>>X^*(P^K). \end{CD}\quad \quad \text{(**)}
\end{equation*}
The horizontal maps in the right-hand square are those defined in A.5 and
A.7, the map $X^{\ast }(S^{K})\rightarrow X^{\ast }(P^{K})$ is that in the
fundamental diagram (A.8), and the map $L(A_{0}\times B_{0})\rightarrow
L(A\times B)$ comes from the inclusion $C(A_{0}\times B_{0})\subset
C(A\times B)$ induced by the reduction map $\End(A\times B)\hookrightarrow %
\End(A_{0}\times B_{0})$. The character $\chi $ of $L(A\times B)$ is that
defined above, and $\chi _{0}$ is the composite of $\chi $ with $%
L(A_{0}\times B_{0})\rightarrow L(A\times B)$.

We begin by drawing some consequences from our condition: 
\begin{equation*}
\Sigma _{0}\cdot D(w_{0})=D(w_{0})\cdot \Sigma _{0}.
\end{equation*}
This condition still holds, even though we are no longer assuming $K$ to be
the Galois closure of $\sigma _{0}E$. In fact, we now assume that $K$ is
large enough to split $\End^{0}(A_{0}\times B_{0})$ (in the sense of A.6).
Note that the condition implies that $\Sigma _{0}\cdot D(w_{0})$ is a
subgroup of $\Sigma _{K}$.

Let $X$ be the set of primes of $K$ dividing $p$. Suppose that the subsets $%
\Sigma _{i}\cdot w_{0}$ and $\Sigma _{j}\cdot w_{0}$ of $X$ have nonempty
intersection. Then $\tau _{i}w_{0}=\tau _{j}w_{0}$ for some $\tau _{i}\in
\Sigma _{i}$ and $\tau _{j}\in \Sigma _{j}$. Hence $\tau _{i}\in \tau
_{j}D(w_{0})$, and so 
\begin{equation*}
\Sigma _{i}\cdot w_{0}=\tau _{i}\Sigma _{0}\cdot w_{0}\subset \tau
_{j}D(w_{0})\Sigma _{0}\cdot w_{0}=\tau _{j}\Sigma _{0}D(w_{0})\cdot
w_{0}=\Sigma _{j}\cdot w_{0}.
\end{equation*}
By symmetry, $\Sigma _{i}\cdot w_{0}\supset \Sigma _{j}\cdot w_{0}$, and so
the two sets are equal: we have shown that the sets $\Sigma _{i}\cdot w_{0}$
and their complex conjugates form a partition of $X.$ Let $X_{0},\ldots
,X_{m-1}$ be the distinct elements of $\{\Sigma _{i}\cdot w_{0}\mid 0\leq
i\leq n-1\}$ with $X_{0}$ chosen to be $\Sigma _{0}w_{0}$, and let 
\begin{equation*}
Y=\{X_{0},\ldots ,X_{m-1},\iota X_{0},\ldots ,\iota X_{m-1}\}.
\end{equation*}
The group $\Sigma _{K}$ acts transitively on $X$ and $Y$, and the
stabilizers of $w_{0}$ and $X_{0}$ are $D(w_{0})$ and $\Sigma _{0}\cdot
D(w_{0})$ respectively. By using $w_{0}$ and $X_{0}$ as base points, we can
identify the map of $\Sigma _{K}$-sets $X\rightarrow Y$ with $\Sigma
_{K}/D(w_{0})\rightarrow \Sigma _{K}/\Sigma _{0}\cdot D(w_{0})$. Each $X_{j}$
then corresponds to the quotient of a left coset of $\Sigma _{0}\cdot
D(w_{0})$ by the right action of $D(w_{0})$. From these remarks, we see that 
\begin{align*}
|X|& =(\Sigma _{K}:D(w_{0})),\quad \\
|Y|& =(\Sigma _{K}:\Sigma _{0}\cdot D(w_{0}))\quad (=2m), \\
|X_{j}|& =(\Sigma _{0}\cdot D(w_{0}):D(w_{0})).
\end{align*}
For $i\in \{0,\ldots ,n-1\}$, define $j(i)$ to be the element of $\{0,\ldots
,m-1\}$ such that $\Sigma _{i}\cdot w_{0}=X_{j(i)}$. For each $j$, there are 
$(\Sigma _{0}\cdot D(w_{0}):\Sigma _{0})=n/m$ sets $\Sigma _{i}$ such that $%
\Sigma _{i}\cdot w_{0}=X_{j}$.

We next compute the terms in the diagram 
\begin{equation*}
\begin{CD} X^{\ast }(L(A))@>>>X^{\ast }(S^{K}) \\ @VVV @VVV \\ X^{\ast
}(L(A_{0}))@>>> X^{\ast }(P^{K}). \end{CD}\quad \quad \text{(**A)}
\end{equation*}
Recall that we have already shown that $X^{\ast }(L(A))\rightarrow X^{\ast
}(S^{K})$ sends the element $[\sigma _{i}]$ of $X^{\ast }(L(A))$ to $\psi
_{i}$.

We use the map $\pi \mapsto f_{\pi }^{K}$ (see A.6, A.7) to identify $%
X^{\ast }(P^{K})$ with 
\begin{equation*}
\{f\in \mathbb{Z}{}^{X}\mid \text{there exists an }m\in \mathbb{Z}{}\text{
such that }f+\iota f=mn_{0}\}.
\end{equation*}
Here $n_{0}=[K_{w_{0}}:\mathbb{Q}{}_{p}]=|D(w_{0})|$. The map $X^{\ast
}(S^{K})\rightarrow X^{\ast }(P^{K})$ is 
\begin{equation*}
f=\sum_{\tau \in \Sigma _{K}}f(\tau )\tau \mapsto \sum_{\tau \in \Sigma
_{K}}f(\tau )\tau w_{0}=\sum_{w\in X}(\sum_{\tau \text{, }\tau
w_{0}=w}f(\tau ))w
\end{equation*}
(see A.8). When $f=\psi _{i}$, $w\in X_{j}$ occurs in the right-hand side
with nonzero coefficient if and only if $j=j(i)$, in which case its
coefficient is $|\Sigma _{0}\cap D(w_{0})|$. Thus the map sends $\psi _{i}$
to $f_{j(i)}$ where $f_{j}$ is the function determined by the conditions 
\begin{eqnarray*}
f_{j}(w) &=&\left\{ 
\begin{array}{ll}
|\Sigma _{0}\cap D(w_{0})|\quad & w\in X_{j} \\ 
0 & w\in X_{j^{\prime }}\text{, }j^{\prime }\neq j
\end{array}
\right. ,\quad \\
f_{j}(w)+f_{j}(\iota w) &=&n_{0},\quad \quad \text{all }w.
\end{eqnarray*}

We identify $X^{\ast }(L(A_{0}))$ with 
\begin{equation*}
\frac{\mathbb{Z}{}^{\Pi _{A_{0}}}}{\{g\mid g=\iota g\text{, }\sum g(\pi )=0\}%
}
\end{equation*}
where $\Pi _{A_{0}}$ is the set of conjugates of $\pi _{A_{0}}$ in $K$ (see
A.7). Let $u=\rho _{0}^{-1}w_{0}$, and let $v_{0}=\sigma _{0}^{-1}w_{0}$.
Note that $\sigma _{i}^{-1}w_{0}$ lies over $u_{0}$ and $(\iota \sigma
_{i})^{-1}w_{0}$ lies over $\iota u_{0}$, $0\leq i\leq n-1$. Using this, we
find that the slope function of the Frobenius germ $\pi _{A_{0}}$ of $A_{0}$
satisfies 
\begin{equation*}
s_{\pi _{A_{0}}}(v)=\left\{ 
\begin{array}{ll}
1/|\Sigma _{E}(v_{0})|\quad & v=v_{0} \\ 
0 & v\text{ lies over }u_{0},\quad v\neq v_{0}
\end{array}
\right.
\end{equation*}
where $\Sigma _{E}(v_{0})=\{\sigma \in \Sigma _{E}\mid \sigma
^{-1}w_{0}=v_{0}\}$ (see A.8). As $s+\iota s=1$, this determines $s$. Note
that 
\begin{equation*}
|\Sigma _{E}(v_{0})|=(\Sigma _{0}\cdot D(w_{0}):\Sigma
_{0})=(D(w_{0}):\Sigma _{0}\cap D(w_{0})).
\end{equation*}
Note also that $X_{0}$ is the set of $w\in X$ lying over the prime $\sigma
_{0}v_{0}$ in $\sigma _{0}E$. For any $\tau \in \Sigma _{i}$ (i.e., such
that $\tau \circ \sigma _{0}=\sigma _{i}$), the diagram 
\begin{equation*}
\begin{array}{ccccccc}
K & 
\begin{smallmatrix}
\tau \\ 
\rightarrow \\ 
\approx
\end{smallmatrix}
& K & \quad & X_{0} & \mapsto & \tau X_{0}=X_{j(i)} \\ 
| &  & | &  &  &  &  \\ 
\sigma _{0}E & 
\begin{smallmatrix}
\quad \\ 
\rightarrow \\ 
\approx
\end{smallmatrix}
& \sigma _{i}E &  & \sigma _{0}v_{0} & \mapsto & \sigma _{i}v_{0}
\end{array}
\end{equation*}
shows that $X_{j(i)}$ is the set of $w\in X$ lying over $\sigma _{i}v_{0}$
in $\sigma _{i}E_{0}$. In other words, $X_{j(i)}$ is the set of $w\in X$
such that $\sigma _{i}^{-1}w=v_{0}$. For $\sigma \in \Sigma _{E}$, $\sigma
\pi _{A_{0}}$ is the Weil germ in $K$ with 
\begin{equation*}
f_{\sigma \pi _{A_{0}}}^{K}(w)=s_{\sigma \pi _{A_{0}}}(w)\cdot n_{0}=s_{\pi
_{A_{0}}}(\sigma ^{-1}w)\cdot n_{0}.
\end{equation*}
When $\sigma =\sigma _{i}$ and $w\in X_{j}$, this becomes 
\begin{equation*}
f_{\sigma _{i}\pi _{A_{0}}}^{K}(w)=\left\{ 
\begin{array}{cc}
|\Sigma _{0}\cap D(w_{0})| & j=j(i) \\ 
0 & j\neq j(i)
\end{array}
.\right.
\end{equation*}
Thus, $f_{\sigma _{i}\pi _{As}}=f_{j(i)}$. In particular, $\sigma _{i}\pi
_{A_{0}}$ depends only on $j(i)$. As the functions $f_{j}$ are distinct, we
see that 
\begin{equation*}
\Pi _{A_{0}}=\{\pi _{0},\ldots ,\pi _{m-1},\iota \pi _{0},\ldots ,\iota \pi
_{m-1}\}
\end{equation*}
where $\pi _{j(i)}=\sigma _{i}\pi _{A_{0}}$. The map $X^{\ast
}(L(A))\rightarrow X^{\ast }(L(A_{0}))$ sends $[\sigma _{i}]$ to $[\pi
_{j(i)}]$, and the map $X^{\ast }(L(A_{0}))\rightarrow X^{\ast }(P^{K})$
sends $[\pi _{j}]$ to $f_{j}$.

We have now computed all the terms in the diagram (**A). It is clear that it
commutes.

We next compute the terms in the diagram 
\begin{equation*}
\begin{CD} X^{\ast }(L(B))@>>> X^{\ast }(S^{K}) \\ @VVV @VVV \\ X^{\ast
}(L(B_{0}))@>>>X^{\ast }(P^{K}). \end{CD}\quad \quad \text{(**B).}
\end{equation*}
Recall that $X^{\ast }(L(B))$ has basis $[\rho _{0}]$, $[\iota \rho _{0}]$,
and that the map $X^{\ast }(L(B))\rightarrow X^{\ast }(S^{K})$ sends $[\rho
_{0}]$ to $\psi $ and $[\iota \rho _{0}]$ to $\iota \psi $. Here $\psi $ is
the characteristic function of $\bigcup \Sigma _{i}$. Clearly, $\mathbb{Q}%
{}[\pi _{B_{0}}]=Q$, and $X^{\ast }(L(B_{0}))=\mathbb{Z}{}^{\Sigma _{Q}}$.
The left-hand vertical map in the diagram is therefore the identity map. Let 
$f\in \mathbb{Z}{}^{X}$ be the function 
\begin{equation*}
f(w)=\left\{ 
\begin{array}{l}
n_{0},\quad \text{if }\rho _{0}^{-1}w=u_{0} \\ 
0,\quad \text{otherwise}
\end{array}
\right. .
\end{equation*}
Then $f\in X(P^{K})$ and the bottom map sends $[\rho _{0}]\mapsto f$ . The
right-hand vertical map sends $\psi $ to $f$.

On combining the diagrams (**A) and (**B), we get the right-hand square in
(**). It remains to compute the kernel of $X^{\ast }(L(A_{0}\times
B_{0}))\rightarrow X^{\ast }(P^{K})$. Note that 
\begin{equation*}
X^{\ast }(L(A_{0}\times B_{0}))=\frac{\mathbb{Z}{}^{\Sigma _{\mathbb{\Pi }%
_{A_{0}}}\sqcup \Sigma _{Q}}}{\{g\mid g=\iota g\text{ and }\sum g(y)=0\}}.
\end{equation*}
The elements $[\pi _{0}],\ldots ,[\pi _{m-1}],[\rho _{0}],[\rho _{0}+\iota
\rho _{0}]$ form a basis for $X^{\ast }(L(A_{0}\times B_{0}))$. They are
mapped respectively to $f_{0},\ldots ,f_{m-1},f,f+\iota f$ in $X^{\ast
}(P^{K}).$ Clearly, 
\begin{equation*}
\frac{n}{m}(f_{0}+\cdots +f_{m-1})+(n-2)f=(n-1)(f+\iota f),
\end{equation*}
and any relation among $f_{0},\ldots ,f_{m-1},f,f+\iota f$ is a multiple of
this one. Therefore, the kernel of $X^{\ast }(L(A_{0}\times
B_{0}))\rightarrow X^{\ast }(P^{K})$ is the free $\mathbb{Z}{}$-module of
rank one generated by 
\begin{equation*}
\chi _{0}=[\frac{n}{m}(\pi _{0}+\cdots +\pi _{m-1})+(n-2)\rho
_{0}-(n-1)(\rho _{0}+\iota \rho _{0})].
\end{equation*}
The map $X^{\ast }(L(A\times B))\rightarrow X^{\ast }(L(A_{0}\times B_{0}))$
sends $\chi $ to $\chi _{0}$, and so we have obtained the diagram (**).

We now prove Theorem 1.5. The group $L(A_{0}\times B_{0})$ acts on the space 
$W(A_{0},B_{0})$ through the characters $\chi _{0}$ and $\iota \chi
_{0}=-\chi _{0}$. Because $\chi _{0}$ is trivial on $P(A_{0}\times B_{0})$, $%
W(A_{0}\times B_{0})$ consists of Tate classes, and because $\chi _{0}$ is
nontrivial on $L(A_{0}\times B_{0})$, the classes are exotic. The group $%
L(A_{0}\times B_{0})$ acts on no other subspace of a space $%
H^{2r}(A_{0}\times B_{0}^{n-2},\mathbb{Q}{}_{\ell }(r))$ through the
character $\chi _{0}$, and so $W(A_{0}\times B_{0})$ contains all the exotic
Tate classes on $A_{0}\times B_{0}^{n-2}$.

From (**), we obtain an exact commutative diagram 
\begin{equation*}
\begin{CD} 0@>>>MT(A\times B)@>>>L(A\times B)@>{\chi}>>T@>>>0\\
@.@AAA@AAA@AA{\cong}A@.\\ 0@>>>P(A_0\times B_0)@>>>L(A_0\times
B_0)@>{\chi_0}>>T_0@>>>0 \end{CD}.
\end{equation*}
It follows that 
\begin{equation*}
P(A_{0}\times B_{0})=L(A_{0}\times B_{0})\cap MT(A\times B)\text{.}
\end{equation*}
If some exotic Hodge class on $A\times B^{n-2}$ is algebraic, then the Hodge
conjecture holds for all powers of $A\times B^{n-2}$ (see Theorem 1.4), and
so (b) of Theorem 1.5 follows from Corollary 1.2.

\begin{remark}
It follows from the above calculations that $P(A_{0})=L(A_{0})$ and $%
P(B_{0})=L(B_{0})$, and so all Tate classes on $A_{0}$ and $B_{0}$ are
Lefschetz.
\end{remark}

\begin{remark}
Choose $E$ to be Galois over $\mathbb{Q}{}$, and identify it with $K$. In
this case, the maps $X^{\ast }(L(A\times B))\rightarrow X^{\ast }(S^{K})$
and $X^{\ast }(L(A_{0}\times B_{0}))\rightarrow X^{\ast }(P^{K})$ are
surjective, and so we obtain an exact commutative diagram 
\begin{equation*}
\begin{CD} 0@>>>\langle\chi\rangle@>>>X^*(L(A\times B))@>>>X^*(S^K)@>>>0\\
@.@VV{\cong}V@VVV@VVV\\ 0@>>>\langle\chi_{0}\rangle@>>>X^*(L(A_0\times
B_0)@>>>X^*(P^K)@>>>0. \end{CD}.
\end{equation*}
The vertical arrows are surjective, and so 
\begin{equation*}
0\rightarrow \Ker(X^{\ast }(L(A\times B))\rightarrow X^{\ast }(L(A_{0}\times
B_{0}))\rightarrow X^{\ast }(S^{K})\rightarrow X^{\ast }(P^{K})\rightarrow 0
\end{equation*}
is exact. Hence 
\begin{equation*}
0\rightarrow P^{K}\rightarrow S^{K}\rightarrow L(A\times B)/L(A_{0}\times
B_{0})
\end{equation*}
is exact, which implies that 
\begin{equation*}
0\rightarrow P^{K}\rightarrow S^{K}\rightarrow L^{K}/T^{K}
\end{equation*}
is exact (notations as in Milne 1999b) because the map $L(A\times
B)/L(A_{0}\times B_{0})\rightarrow L^{K}/T^{K}$ is injective. Therefore 
\begin{equation*}
P^{K}=S^{K}\cap T^{K}\quad \text{(intersection inside\textrm{\ }}L^{K}),
\end{equation*}
and we recover ibid., Theorem 6.1.
\end{remark}

\begin{appendix}%
%

\section{Abelian Varieties with Many Endomorphisms}

\subsection{A.1. Notations}

Throughout, $\mathbb{Q}^{\al}$ is the algebraic closure of $\mathbb{Q}$ in $%
\mathbb{C}$, and $\Gamma =\Gal(\mathbb{Q}^{\al}/\mathbb{Q})$. Complex
conjugation on $\mathbb{C}$, or a subfield of $\mathbb{C}$, is denoted by $%
\iota $ or $x\mapsto \bar{x}$. In A.8, we fix a prime $w_{0}$ of $\mathbb{Q}%
{}^{\text{al}}$ dividing $p$, and denote the residue field at $w_{0}$ by $%
\mathbb{F}{}$. We denote the restriction of $w_{0}$ to a subfield of $%
\mathbb{Q}^{\text{al}}$ by the same symbol. For a finite \'{e}tale $\mathbb{Q%
}$-algebra $E$, $\Sigma _{E}=\Hom(E,\mathbb{Q}^{\al})$. For a subfield $K$
of $\mathbb{Q}{}^{\text{al}}$ Galois over $\mathbb{Q}{}$, $\Sigma _{K}$ can
be identified with $\Gal(K/\mathbb{Q}{})$.

A \emph{CM-algebra} $E$ is a finite product of finite field extensions of $%
\mathbb{Q}$ admitting an involution $\iota _{E}$ that is nontrivial on each
factor and such that $\sigma (\iota _{E}x)=\overline{\sigma (x)}$ for all $%
\sigma \colon E\rightarrow \mathbb{C}$; equivalently, $E$ is a finite
product of CM-fields.

For a finite set $Y$, $\mathbb{Z}^{Y}$ denotes the set of functions $f\colon
Y\rightarrow \mathbb{Z}$. We sometimes denote such a function by $\sum f(y)y$%
.

For a group of multiplicative type $T$ over $\mathbb{Q}$, $X^{\ast }(T)%
\overset{\text{df}}{=}\Hom(T_{\mathbb{Q}^{\al}},\mathbb{G}_{m})$ is the
character group. We often use the pairing 
\begin{equation*}
\chi ,\mu \mapsto \langle \chi ,\mu \rangle \overset{\text{df}}{=}\chi \circ
\mu \colon X^{\ast }(T)\times X_{\ast }(T)\rightarrow \End(\mathbb{G}%
_{m})\cong \mathbb{Z}
\end{equation*}
to identify the cocharacter group $X_{\ast }(T)\overset{\text{df}}{=}\func{%
Hom}_{\mathbb{Q}{}^{\text{al}}}(\mathbb{\mathbb{G}{}}_{m},T_{\mathbb{Q}{}^{%
\text{al}}})$ of $T$ with the $\mathbb{Z}$-linear dual of $X^{\ast }(T)$.

Let $\rho \colon T\rightarrow \GL(V)$ be a representation of a group $T$ of
multiplicative type on a finite-dimensional $\mathbb{Q}$-vector space $V$.
For any subfield $\Omega $ of $\mathbb{C}$ that splits $T$, there is a
decomposition 
\begin{equation*}
V\otimes _{\mathbb{Q}}\Omega \cong \bigoplus_{\chi \in X^{\ast }(T)}V_{\chi }
\end{equation*}
where $V_{\chi }$ is the subspace of $V\otimes _{\mathbb{Q}}\Omega $ on
which $T$ acts through $\chi $. If $V_{\chi }$ is nonzero, then we say that $%
\chi $ \emph{occurs }in $V$. When $\Omega $ is Galois over $\mathbb{Q}$, a
subspace $\bigoplus_{\chi \in \Xi }V_{\chi }$, $\Xi \subset X^{\ast }(T)$,
is defined over $\mathbb{Q}$ (i.e., of the form $W\otimes _{\mathbb{Q}%
{}}\Omega $ for some subspace $W\subset V$) if and only if $\Xi $ is stable
under $\Gamma $. The subspace of vectors in $V$ fixed by $T$ (in the sense
of Milne 1999a, \S 3) is denoted $V^{T}$.

For a finite \'{e}tale $\mathbb{Q}$-algebra $E$, $(\mathbb{G}_{m})_{E/%
\mathbb{Q}}=_{\text{df}}\Res_{E/\mathbb{Q}}(\mathbb{G}_{m})$ (Weil
restriction of scalars), so that $X^{\ast }((\mathbb{G}_{m})_{E/\mathbb{Q}})=%
\mathbb{Z}^{\Sigma _{E}}$. Under this identification, an element $f=\sum
f(\sigma )\sigma $ of $\mathbb{\mathbb{Z}{}}^{\Sigma _{E}}$ maps an element $%
a$ of $E^{\times }=(\mathbb{\mathbb{G}{}}_{m})_{E/\mathbb{Q}{}}(\mathbb{Q}{})
$ to $a^{f}=\prod (\sigma a)^{f(\sigma )}$. We sometimes identify a subset $%
\Delta $ of $\Sigma _{E}$ with the character $\sum_{\sigma \in \Delta
}\sigma $; for example, if $V$ is an $E$-vector space, then $(V^{\otimes
r}\otimes \Omega )_{\Delta }$ is the subspace on which $a\in E$ acts as $%
\prod_{\sigma \in \Delta }\sigma a$.

There is a natural correspondence\footnote{%
Experts will recognize the Tannakian significance of this correspondence
(Saavedra 1972, V 3.1.4; Deligne and Milne 1982, p.\ 190).} between 
\begin{mylist}
\item triples $(T,w,t)$ comprising a group of multiplicative type $
T$ over $\mathbb{Q}$, a 
cocharacter $w$ of $T$, and a character $t$ such that $t\circ w=-2$; and 
\item pairs $(T_0,\varepsilon )$ comprising a group of multiplicative type $
T_0$ and an 
element $\varepsilon$ of order $1$ or $2$ in $T_0(\mathbb{Q})$.  
\end{mylist}
Given $(T,w,t)$, define $T_{0}$ to be the kernel of $t$ and $\varepsilon $
to be $w(-1)$. Conversely, given $(T_{0},\varepsilon )$, define $T$ by the
diagram 
\begin{equation*}
\begin{CD} 0@>>>\mu_2@>>>{\Bbb{G}}_m@>{-2}>>{\Bbb{G}}_m@>>>0\\
@.@VV{\varepsilon}V@VV{w}V@VV{=}V\\ 0@>>>T_0@>>>T@>{t}>>{\Bbb{G}}_m@>>>0
\end{CD}
\end{equation*}
in which $T=(T_{0}\times \mathbb{G}_{m})/\mu _{2}$. If $(T_{1},w_{1},t_{1})%
\subset (T_{2},w_{2},t_{2})$, then $T_{1}=T_{2}$ if and only if $%
(T_{1})_{0}=(T_{2})_{0}$.

Let $\rho _{0}\colon T_{0}\rightarrow GL(V)$ be a representation of $T_{0}$
such that $\rho _{0}(\epsilon )$ acts on $V$ as multiplication by the scalar 
$-1$, and let $W$ be a one-dimensional vector space with basis $e$. Then 
\begin{equation*}
(x,y)\mapsto (\rho _{0}(x)\cdot y,y^{-2})\colon T_{0}\times \mathbb{G}%
{}_{m}\rightarrow \GL(V)\times \GL(W)
\end{equation*}
sends $(\epsilon ,\epsilon )$ to $1$, and therefore defines a homomorphism $%
\rho \colon T\rightarrow \GL(V)\times \GL(W)$. Note that $(\rho \circ w)(y)$
acts on $V$ as $y$, and that the composite of $\rho $ with the projection to 
$\GL(W)$ is $t.$ Let $s\in V^{\otimes i}$. If $s$ is fixed by $T$, then $i$
is even. There is a one-to-one correspondence 
\begin{equation*}
s\leftrightarrow s\otimes e^{\otimes j}
\end{equation*}
between the elements $s$ of $V^{\otimes 2j}$ fixed by $T_{0}$ and the
elements of $V^{\otimes 2j}\otimes W^{\otimes j}$ fixed by $T$.

For a smooth projective variety $X$, $\mathcal{Z}^{r}(X)$ is the space of
algebraic cycles on $X$ of codimension $r$ with coefficients in $\mathbb{Q}$%
, and $\mathcal{Z}_{\text{num}}^{r}(X)$ is the quotient of $\mathcal{Z}%
^{r}(X)$ by numerical equivalence. The space $\mathcal{Z}_{\text{num}}^{\ast
}(X)=\bigoplus_{r}\mathcal{Z}_{\text{num}}^{r}(X)$ becomes a $\mathbb{Q}$%
-algebra under intersection product. An \emph{algebraic class} in a
cohomology group with coefficients in a field $\Omega $ is an element of the 
$\Omega $-subspace spanned by the classes of algebraic cycles.

For an abelian variety $A$ over an algebraically closed field $k$ of
characteristic zero, we often implicitly assume that there is given an
embedding $\sigma \colon k\rightarrow \mathbb{C}{}$ so that we can define $%
H^{r}(A,\mathbb{Q}{})$ to be $r$th cohomology group of the complex manifold $%
(\sigma A)(\mathbb{C}{})$. We let $\func{Hom}^{0}(A,B)=\func{Hom}%
(A,B)\otimes _{\mathbb{Z}{}}\mathbb{Q}{}$.

For Hodge structures and class field theory, we follow the usual conventions
of those areas rather than the conventions of Deligne used in my previous
papers. For example, $z\in \mathbb{C}{}^{\times }$ acts on a Hodge structure
of type $(r,s)$ as $z^{r}\bar{z}^{s}$, and the Artin reciprocity maps send
prime elements to the Frobenius element $x\mapsto x^{q}$.

We sometimes use $[x]$ to denote an equivalence class containing $x$, and $%
|X|$ to denote the order of a finite set $X$.

For an explanation of the various cohomology groups of varieties, and their
Tate twists, see Deligne 1982, \S 1.

This section summarizes results due to many mathematicians. Omitted proofs
can be found in Milne 1999a, 1999b, Tate 1968/69, or in the references for
those articles.

\subsection{A.2. Generalities}

Let $A$ be an abelian variety over an algebraically closed field $k$. The
reduced degree\footnote{%
Let $R$ be a semisimple algebra of finite degree over $\mathbb{Q}$. Then $R$
is a product of simple algebras, say, $R=R_{1}\times \cdots \times R_{m}$,
and the centre $E_{i}$ of each $R_{i}$ is a field. The \emph{reduced } \emph{%
degree} $[R:\mathbb{Q}]_{\text{red}}$ of $R$ over $\mathbb{Q}$ is defined to
be $\sum_{i=1}^{m}[R_{i}:E_{i}]^{\frac{1}{2}}[E_{i}:\mathbb{Q}]$.} of the $%
\mathbb{Q}$-algebra $\End^{0}(A)$ is $\leq 2\dim A$, and when equality holds
the abelian variety is said\footnote{%
Often such an abelian variety is said to admit ``complex multiplication'',
but this conflicts with classical terminology --- see Lange and Birkenhake
1992, p.\ 268. Also ``multiplication'' for ``endomorphism'' seems archaic.}
to have \emph{many endomorphisms\/}. An isotypic\footnote{%
An abelian variety is said to be $\emph{isotypic}$ if it is isogenous to a
power of a simple abelian variety.} abelian variety has many endomorphisms
if and only if $\End^{0}(A)$ contains a field of degree $2\dim A$ over $%
\mathbb{Q}{}$, and an arbitrary abelian variety has many endomorphisms if
and only each isotypic isogeny factor of it does. Equivalent conditions:

\begin{enumerate}
\item  \ the $\mathbb{Q}$-algebra $\End^{0}(A)$ contains an \'{e}tale
subalgebra of degree $2\dim A$ over $\mathbb{Q}{}$;

\item  \ for a Weil cohomology $X\mapsto H^{\ast }(X)$ with coefficient
field $\Omega $, the centralizer of $\End^{0}(A)$ in $\End_{\Omega
}(H^{1}(A))$ is commutative (in which case it equals $C(A)\otimes _{\mathbb{Q%
}{}}\Omega $ where $C(A)$ is the centre of $\End^{0}(A)$);

\item  \ (characteristic zero) $A$ has CM-type, i.e., its Mumford-Tate group
(see A.3 below) is commutative (hence a torus);

\item  \ (characteristic $p\neq 0$) $A$ is isogenous to an abelian variety
defined over $\mathbb{F}$ (theorems of Tate and Grothendieck).
\end{enumerate}

\ Let $k\subset k^{\prime }$ be algebraically closed fields. The functor $%
A\mapsto A_{k^{\prime }}$ from the category of abelian varieties over $k$ to
the similar category over $k^{\prime }$ is fully faithful, because the map
on torsion points $A(k)_{\text{tors}}\rightarrow A(k^{\prime })_{\text{tors}%
} $ is bijective and $A(k)_{\text{tors}}$ is Zariski dense in $A$. That the
functor becomes essentially surjective on the categories of abelian
varieties with many endomorphisms up to isogeny is a result of Grothendieck
(Oort 1973). Thus, in large part, the theory of abelian varieties with many
endomorphisms up to isogeny over an algebraically closed field depends only
on the characteristic of the field.

\subsection{A.3. The groups attached to an abelian variety with many
endomorphisms}

Let $A$ be an abelian variety with many endomorphisms over an algebraically
closed field $k$, and let $C(A)$ be the centre of $\End^{0}(A)$. Every
Rosati involution on $\End^{0}(A)$ stabilizes $C(A)$, and the different
Rosati involutions restrict to the same involution on $C(A)$, which we
denote $^{\dagger }$. Each factor of $C(A)$ is either a CM-field, on which $%
^{\dagger }$ acts as complex conjugation, or is $\mathbb{Q}{}$.

\subsubsection{The Lefschetz group}

We define $L(A)_{0}$ to be the group of multiplicative type over $\mathbb{Q}$
such that, for all commutative $\mathbb{Q}$-algebras $R$, 
\begin{equation*}
L(A)_{0}(R)=\{\alpha \in C(A)\otimes R\mid \alpha \alpha ^{\dagger }=1\}.
\end{equation*}
Let $\varepsilon =-1\in L(A)_{0}(\mathbb{Q})$, and let $(L(A),w,t)$ be the
triple associated (as in A.1) with $(L(A)_{0},\varepsilon )$.

Then 
\begin{equation*}
L(A)(\mathbb{Q}{})\cong \{\alpha \in C(A)^{\times }\mid \alpha \alpha
^{\dagger }\in \mathbb{Q}{}^{\times }\},
\end{equation*}
and, on $\mathbb{Q}{}$-points, $w$ is $x\mapsto x$ and $t$ is $x\mapsto
(xx^{\dagger })^{-1}$.

\subsubsection{The motivic group}

Because $L(A)_{0}$ is a subgroup of $\End^{0}(A)^{\times }$, it acts on $%
\mathcal{Z}_{\text{num}}^{\ast }(A^{s})$ for all $s$, and we define $%
M(A)_{0} $ to be the largest algebraic subgroup of $L(A)_{0}$ acting
trivially on these $\mathbb{Q}$-algebras. Then $-1\in M(A)_{0}(\mathbb{Q}{})$%
, and we let $(M(A),w,t)$ be the triple associated with $(M(A)_{0},-1)$.

\subsubsection{The Mumford-Tate group}

When $k$ has characteristic zero, $L(A)_{0}$ acts on the $\mathbb{Q}$%
-algebra of Hodge classes on $A^{s}$ for all $s$, and we define $MT(A)_{0}$
to be the subgroup of $L(A)_{0}$ fixing the elements of these $\mathbb{Q}$%
-algebras. Again $-1\in MT(A)_{0}$, and we let $(MT(A),w,t)$ be the triple
associated with $(MT(A)_{0},-1)$.

\subsubsection{The group $P$}

Let $k=\mathbb{F}$. A model $A_{1}$ of $A$ over a finite field $\mathbb{F}%
_{q}$ defines a Weil $q$-number $\pi _{1}$, whose class $\pi _{A}$ in $%
W(p^{\infty })$ (see A.6 below) is independent of the choice of $A_{1}$. The
group $P(A)$ is defined to be the smallest algebraic subgroup of $L(A)$
containing some power of $\pi _{1}$ --- again, it is independent of the
choice of $A_{1}$.

Let $\pi _{1}$ be a Weil $p^{2n}$-number representing $\pi _{A}$. Then $\pi
_{1}/p^{n}\in L(A)_{0}$, and $P(A)_{0}$ is the smallest algebraic subgroup
of $L(A)_{0}$ containing some power of $\pi _{1}/p^{n}$.\smallskip

Let $H^{\ast }$ be a Weil cohomology with coefficients in a field $\Omega $.
Since $L(A)_{0}\subset (\mathbb{\mathbb{G}{}}_{m})_{E/\mathbb{Q}{}}$, there
is a natural action of $L(A)_{0}$ on $H^{1}(A,\mathbb{\Omega })$, and $%
\varepsilon $ acts as $-1$. Hence (see A.1) there is a natural action of $%
L(A)$ on 
\begin{equation*}
H^{r}(A^{s},\mathbb{\Omega }{})(m)\cong (\bigwedge^{r}(\bigoplus_{s\text{
copies}}H^{1}(A,\mathbb{\Omega }{})))\otimes (\Omega {}(1))^{\otimes m}\text{%
.}
\end{equation*}

\begin{lemma*}
Let $A$ be an abelian variety with many endomorphisms over an algebraically
closed field $k$, and let $H^{\ast }$ be a Weil cohomology with coefficients
in a field $\Omega $. Let $H^{2\ast }(A^{s})(\ast
)=\bigoplus_{r}H^{2r}(A^{s})(r)$. Then, for all $s$,

\begin{enumerate}
\item  \ $H^{2\ast }(A^{s})(\ast )^{L(A)}$ is the $\Omega $-subalgebra of $%
H^{2\ast }(A^{s})(\ast )$ generated by the classes of divisors on $A^{s}$
(i.e., it is the space of Lefschetz classes);

\item  \ $H^{2\ast }(A^{s})(\ast )^{M(A)}$ is the space of algebraic classes
in $H^{2\ast }(A^{s})(\ast )$, provided numerical equivalence coincides with
homological equivalence for $H^{\ast }$;

\item  $H^{2\ast }(A^{s})(\ast )^{MT(A)}$ is the space of Hodge classes on $%
A^{s}$ when $k$ has characteristic zero and $H^{\ast }$ is the cohomology
defined by an embedding $k\rightarrow \mathbb{C}{}$;

\item  $H^{2\ast }(A^{s},\mathbb{Q}{}_{\ell }(\ast ))^{P(A)}$ is the space
of $\ell $-adic Tate classes on $A^{s}$ when $k=\mathbb{F}{}$ and $H^{\ast }$
is $\ell $-adic \'{e}tale cohomology.
\end{enumerate}
\end{lemma*}

\begin{proof}
Statement (a) is proved in Milne 1999a (Theorem 4.4).

For (b), recall that theorems of Jannsen and Deligne show that the category
of abelian motives over $k$, defined using the numerical equivalence classes
of algebraic cycles as correspondences, is Tannakian (Jannsen 1992). Almost
by definition, $M(A)$ is the fundamental group of the Tannakian subcategory $%
\langle A\rangle ^{\otimes }$ of this category generated by $A$ and the Tate
object. When numerical equivalence coincides with homological equivalence,
there is a natural map 
\begin{equation*}
\Hom(\1,h^{2r}(A^{s})(r))\otimes _{\mathbb{Q}}\Omega \rightarrow \Hom%
_{\Omega }(\Omega ,H^{2r}(A^{s})(r))^{M(A)},
\end{equation*}
which the theory of Tannakian categories shows to be bijective. But $\Hom(\1%
,h^{2r}(A^{s})(r))=\mathcal{Z}_{\text{num}}^{r}(A^{s})$.

Statement (c) is proved in Deligne 1982 (see the proof of 3.4).

Almost by definition of $P(A)$, $H^{2\ast }(A^{s},\mathbb{Q}{}_{\ell }(\ast
))^{P(A)}$ consists of the classes fixed by the Frobenius germ $\pi _{A}$,
and these are exactly the Tate classes.
\end{proof}

Thus, under the hypotheses in each part of the lemma, knowing the group $%
?(A)_{\Omega }$ is equivalent to knowing the corresponding spaces of fixed
classes: $?(A)_{\Omega }$ is the largest algebraic subgroup of $\GL%
(H^{1}(A))\times \mathbb{G}_{m}$ fixing the particular classes on all $A^{s}$%
, and the particular classes are exactly those fixed by $?(A)_{\Omega }$.

\begin{theorem*}
\begin{enumerate}
\item  For any abelian variety $A$ with many endomorphisms over an
algebraically closed field $k$ of characteristic zero, $MT(A)\subset
M(A)\subset L(A)$, and

\begin{enumerate}
\item  the Hodge conjecture holds for all powers of $A$ if and only if $%
MT(A)=M(A)$;

\item  all Hodge classes on all powers of $A$ are Lefschetz if and only if $%
MT(A)=L(A)$.
\end{enumerate}

When $k=\mathbb{Q}^{\al}$, ``Hodge'' can be replaced by ``Tate'' in the
above statements.

\item  For any abelian variety $A_{0}$ over $\mathbb{F}$, $P(A_{0})\subset
M(A_{0})\subset L(A_{0})$, and

\begin{enumerate}
\item  all $\ell $-adic Tate classes on all powers of $A_{0}$ are algebraic
for one (or all) $\ell $ if and only if $P(A_{0})=M(A_{0})$;

\item  all $\ell $-adic Tate classes on all powers of $A_{0}$ are Lefschetz
for one (or all) $\ell $ if and only if $P(A_{0})=L(A_{0})$.
\end{enumerate}
\end{enumerate}
\end{theorem*}

\begin{proof}
Since every character of $L(A)$ occurs in a space of the form $%
H^{r}(A^{s})(m)$, we see that the subgroups of $L(A)$ are determined by
their invariants in these spaces. Thus (a) of the theorem is an immediate
consequence of the lemma. That ``Hodge'' can be replaced by ``Tate'' follows
from Pohlmann 1968.

If $P(A_{0})=M(A_{0})$, then the lemma shows that the $\ell $-adic Tate
conjecture holds for all powers of $A_{0}$ and all $\ell $ in the set in
Proposition B.2, but if the $\ell $-adic Tate conjecture holds for one $\ell 
$ then it holds for all (Tate 1994, 2.9). Conversely, if the $\ell $-adic
Tate conjecture holds for all powers of $A_{0}$ and a single $\ell $, then
numerical equivalence coincides with $\ell $-homological equivalence for
that $\ell $ (Milne 1986, 8.4), and the preceding lemma then shows that $%
P(A_{0})_{\mathbb{Q}{}_{\ell }}=M(A_{0})_{\mathbb{Q}{}_{\ell }}$. As $%
P(A_{0})\subset M(A_{0})$, this implies that $P(A_{0})=M(A_{0})$.

The proof of the remaining statement is similar.
\end{proof}

\begin{example*}
If $A$ has dimension $1$, then either $\End^{0}(A)$ is a quadratic imaginary
field $E$ or a quaternion algebra $D$ with centre $\mathbb{Q}{}$. In the
first case, all the groups attached to $A$ equal $(\mathbb{\mathbb{G}{}}%
_{m})_{E/\mathbb{Q}{}}$ and in the second, all the groups attached to $A$
equal $\mathbb{\mathbb{G}{}}_{m}.$ Hence, there are no exotic Hodge or Tate
classes on any power of an elliptic curve, and the Hodge and Tate
conjectures hold.
\end{example*}

\subsection{A.4. Classification over $\mathbb{C}$ of abelian varieties with
many endomorphisms}

Let $E$ be a CM-algebra. A \emph{CM-type\/} on $E$ is the choice of one out
of every pair of complex conjugate homomorphisms $E\rightarrow \mathbb{C}{}$%
. It can variously be considered as:

\begin{enumerate}
\item  \ a partition $\Sigma _{E}=\Phi \cup \iota \Phi $;

\item  \ a function $\varphi \colon \Sigma_E\to \mathbb{Z}$ such that, for
all $\sigma$, $\varphi (\sigma )\geq 0$ and $\varphi (\sigma )+\varphi
(\iota\sigma )=1$;

\item  \ the choice of an isomorphism $E\otimes _{\mathbb{Q}}\mathbb{R}%
\rightarrow \mathbb{\mathbb{C}{}}^{\Sigma _{F}}$ where $F$ is the product of
the maximal real subfields of the factors of $E$.
\end{enumerate}

\noindent Here $\Phi $ is the support of $\varphi $ and $\varphi $ is the
characteristic function of $\Phi $.

Let $A$ be a simple abelian variety over $\mathbb{C}$ with many
endomorphisms. Then $\End^{0}(A)$ is a CM-field $E$, and the action of $E$
on $\Gamma (A,\Omega ^{1})$ defines a CM-type $\Phi $ on $E$, which is
primitive, i.e., not the extension of a CM-type on a proper CM-subfield of $%
E $. The map $A\mapsto (E,\Phi )$ defines a bijection from the set of
isogeny classes of simple abelian varieties over $\mathbb{C}$ with many
endomorphisms to the set of isomorphism classes of pairs $(E,\Phi )$. It
remains to classify the pairs $(E,\Phi )$.

Fix a (large) CM-field $K\subset \mathbb{Q}^{\al}$, finite and Galois over $%
\mathbb{Q}$. The \emph{Serre group } $S^{K}$ of $K$ is the quotient of $(%
\mathbb{G}_{m})_{K/\mathbb{Q}}$ whose character group consists of the $f\in 
\mathbb{Z}{}^{\Sigma _{K}}$ for which there is an integer $wt(f)$ (the \emph{%
weight }of $f)$ such that $f(\tau )+f(\iota \tau )=wt(f)$ for all $\tau \in
\Sigma _{K}$, that is, 
\begin{equation*}
X^{\ast }(S^{K})=\{f\in \mathbb{Z}^{\Sigma _{K}}\mid f+\iota f\text{ is
constant}\}.
\end{equation*}

The \emph{reflex field }of $(E,\Phi )$ is the fixed field of the subgroup $%
\{\tau \in \Gamma \mid \tau \Phi =\Phi \}$ of $\Gamma $. We classify the
pairs $(E,\Phi )$ whose reflex field is contained in $K$. Let $\varphi $ be
the characteristic function of $\Phi $. For each $\sigma \colon E\rightarrow 
\mathbb{Q}^{\al}$ and $\tau \in \Gal(\mathbb{Q}^{\al}/\mathbb{Q})$, define 
\begin{equation*}
\psi _{\sigma }(\tau )=\varphi (\tau ^{-1}\circ \sigma ).
\end{equation*}
Then $\psi _{\sigma }(\tau )$ depends only on $\tau |K$, and for any $\rho
\in \Gal(\mathbb{Q}^{\al}/\mathbb{Q})$, $\psi _{\rho \circ \sigma }=\rho
\psi _{\sigma }$. Thus, $\{\psi _{\sigma }\}$ is a $\Gamma $-orbit in $%
\mathbb{Z}{}^{\Sigma _{K}}$. The map $(E,\Phi )\mapsto \{\psi _{\sigma }\}$
is a bijection from the set of isomorphism classes of pairs $(E,\Phi )$
comprising a CM-field and a primitive CM-type whose reflex field is
contained in $K$ to the set of $\Gamma $-orbits of elements $f$ of $X^{\ast
}(S^{K})$ such that $f(\tau )\geq 0$ for all $\tau $ and $wt(f)=1$.

\subsection{A.5. Calculation of the groups over $\mathbb{C}$}

Let $A$ be an abelian variety with many endomorphisms over $\mathbb{C}$.
Then $A$ is isogenous to a product $A_{1}^{s_{1}}\times \cdots \times
A_{t}^{s_{t}}$ with the $A_{i}$ simple and pairwise nonisogenous, and 
\begin{eqnarray*}
L(A) &\cong &L(A_{1}\times \cdots \times A_{t}),\text{ (in fact }%
L(A)_{0}\cong L(A_{1})_{0}\times \cdots \times L(A_{t})_{0}\text{)} \\
M(A) &\cong &M(A_{1}\times \cdots \times A_{t}) \\
MT(A) &\cong &MT(A_{1}\times \cdots \times A_{t}).
\end{eqnarray*}
Thus, in the following, we assume that $A$ is a product of pairwise
nonisogenous simple abelian varieties. Then, $E=_{\text{df}}\End^{0}(A)$ is
a CM-algebra. The action of $E$ on $H^{1,0}(A)$ defines a CM-type $\Phi $ on 
$E$, and the Rosati involution is $\iota _{E}$.

\subsubsection{The Lefschetz group}

The group $L(A)$ is the subgroup of $(\mathbb{G}_{m})_{E/\mathbb{Q}}$ whose
character group is 
\begin{equation*}
\frac{\mathbb{Z}^{\Sigma _{E}}}{\{g\mid g=\iota g\text{ and }\sum g(\sigma
)=0\}}.
\end{equation*}
The weight map $w\colon \mathbb{\mathbb{G}{}}_{m}\rightarrow L(A)$
corresponds to the map 
\begin{equation*}
\lbrack g]\mapsto wt(g)\overset{\text{df}}{=}\sum_{\sigma \in \Sigma
_{E}}g(\sigma )
\end{equation*}
on characters, and the homomorphism $t\colon L(A)\rightarrow \mathbb{\mathbb{%
G}{}}_{m}$ giving the action of $L(A)$ on the Tate object $\mathbb{Q}{}(1)$
sends $1\in X^{\ast }(\mathbb{G}{}_{m})$ to the element of $X^{\ast }(L(A))$
represented by $-\sigma -\iota \sigma $ for any $\sigma \in \Sigma _{E}$.

The group $L(A)_{0}$ is the subgroup of $(\mathbb{G}_{m})_{E/\mathbb{Q}}$
whose character group is 
\begin{equation*}
\frac{\mathbb{Z}^{\Sigma _{E}}}{\{g\mid g=\iota g\}}.
\end{equation*}
The map $\mu _{2}\rightarrow L(A)_{0}$ corresponds to the map on characters $%
[g]\mapsto \sum g(\sigma )\mod2.$

When $A$ is simple, the map $\sigma \mapsto \psi _{\sigma }$ is bijective
and commutes with the action of $\Gamma $, and so it identifies $L(A)$ with
the torus whose character group is 
\begin{equation*}
\frac{\mathbb{Z}^{\Psi }}{\{g\mid g=\iota g\text{ and }\sum g(\psi )=0\}}%
,\quad \Psi =\{\psi _{\sigma }\mid \sigma \in \Sigma _{E}\}.
\end{equation*}

\subsubsection{The Mumford-Tate group}

The Hodge decomposition on $H^{r}(A^{s},\mathbb{Q}{})(m)$ is defined over $%
\mathbb{Q}^{\al}$, i.e., there is a decomposition 
\begin{equation*}
H^{r}(A^{s},\mathbb{Q}{})(m)\otimes \mathbb{Q}{}^{\text{al}}\cong
\bigoplus_{i+j=r-2m}H^{r}(A^{s})(m)^{i,j}
\end{equation*}
that becomes the Hodge decomposition when tensored with $\mathbb{C}{}$.
Since $L(A)_{0}\subset (\mathbb{\mathbb{G}{}}_{m})_{E/\mathbb{Q}{}}$, there
is a natural action of $L(A)_{0}$ on $H^{1}(A,\mathbb{Q}{})$, and $%
\varepsilon $ acts as $-1$. Hence (see A.1) there is a natural action of $%
L(A)$ on 
\begin{equation*}
H^{r}(A^{s},\mathbb{Q}{})(m)\cong (\bigwedge^{r}(\bigoplus_{s\text{ copies}%
}H^{1}(A,\mathbb{Q}{})))\otimes (\mathbb{Q}{}(1))^{\otimes m}\text{.}
\end{equation*}
For $\chi =[g]\in X^{\ast }(L(A))$, $(H^{r}(A^{s})(m))_{\chi }$ is of Hodge
type 
\begin{equation*}
(\sum_{\sigma \in \Phi }g(\sigma ),\sum_{\sigma \in \iota \Phi }g(\sigma )).
\end{equation*}
Every character of $L(A)$ occurs in $H^{r}(A^{s},\mathbb{Q}{})(m)$ for some $%
r,s,m$, and if $[g]$ occurs in $H^{r}(A^{s})(m)$, then $wt(g)=r-2m$. A
character $\chi $ of $L(A)$ is trivial on $MT(A)$ if and only if $\oplus
_{\tau \in \Gamma }H^{2r}(A^{s})(r)_{\tau \chi }$ is purely of type $(0,0)$
for some $r,s$ for which the space is nonzero. Hence, a character $\chi =[g]$
of $L(A)$ is trivial on $MT(A)$ if and only if 
\begin{equation*}
\sum_{\sigma \in \Phi }g(\tau \circ \sigma )=0\quad \text{for all\textrm{\ }}%
\tau \in \Gamma .
\end{equation*}

\subsubsection{The motivic group}

Let $\chi \in X^{\ast }(L(A))$. Then $\chi $ is trivial on $M(A)$ if and
only if $H^{2r}(A^{s})(r)_{\chi }$ contains a nonzero algebraic class for
some $r$ and $s$, in which case all the spaces $H^{2r}(A^{s})(r)_{\chi }$
consist entirely of algebraic classes (see (b) of the lemma in A.3).

\subsubsection{Second description of $MT(A)$}

There is another description of $MT(A)$ that is useful. Let $K$ be a
CM-subfield of $\mathbb{Q}{}^{\text{al}}$, finite and Galois over $\mathbb{Q}%
{}$, and let $S^{K}$ be its Serre group. Let $\tau _{0}\in \Sigma _{K}$ be
the given embedding of $K$ into $\mathbb{Q}^{\al}$. Then $f\mapsto f(\tau
_{0})$ is a cocharacter $\mu ^{K}$ of $S^{K}$ with the property that $\mu
^{K}+\iota \mu ^{K}$ is fixed by $\Gamma $ and so is defined over $\mathbb{Q}%
{}$. The pair $(S^{K},\mu ^{K})$ is universal: if $T$ is a second torus over 
$\mathbb{Q}$ and $\mu \in X_{\ast }(T)$ is defined over $K$ and $\mu +\iota
\mu $ is defined over $\mathbb{Q}$, then there is a unique homomorphism $%
\rho _{\mu }\colon S^{K}\rightarrow T$ such that $(\rho _{\mu })_{\mathbb{Q}%
^{\al}}\circ \mu ^{K}=\mu $. On characters, $\rho _{\mu }$ sends $\chi \in
X^{\ast }(T)$ to the element $f$ of $X^{\ast }(S^{K})$ with $f(\tau )=$ $%
\langle \chi ,\tau \mu \rangle $ for all $\tau $.

Let $A$ be an abelian variety of CM-type $(E,\Phi )$, and let $\mu _{\Phi }$
be the cocharacter of $L(A)$ sending a character $[g]$ of $L(A)$ to $%
\sum_{\sigma \in \Phi }g(\sigma )$. If $K$ contains the reflex field of $%
\Phi $, then $\mu _{\Phi }$ is defined over $K$. Moreover $\mu _{\Phi
}+\iota \mu _{\Phi }$ is $[g]\mapsto wt(g)$, which is defined over $\mathbb{Q%
}$, and so there is a unique homomorphism $\rho _{\Phi }\colon
S^{K}\rightarrow L(A)$ such that $\rho _{\Phi }\circ \mu ^{K}=\mu _{\Phi }.$
It sends a character $g$ of $L(A)$ to the character $f$ of $S^{K}$ such that 
\begin{equation*}
f(\tau )=\langle \lbrack g],\tau \mu _{\Phi }\rangle =\langle \tau
^{-1}[g],\mu _{\Phi }\rangle =\sum_{\sigma \in \Phi }g(\tau \circ \sigma ).
\end{equation*}
The image of this homomorphism is $MT(A)$. It is obvious that this
description agrees with the previous one.

\subsection{A.6. Classification over $\mathbb{F}$ of abelian varieties}

A \emph{Weil} $q$\emph{-number of weight} $m$ is an element $\pi $ of a
field of characteristic zero such that $q^{N}\pi $ is an algebraic integer
for some $N$ and $\sigma (\pi )\cdot \iota (\sigma (\pi ))=q^{m}$ for all
homomorphisms $\sigma \colon \mathbb{Q}[\pi ]\hookrightarrow \mathbb{C}$.
The conditions imply that $q^{N}\pi $ is a unit at all finite primes $v$ of $%
\mathbb{Q}{}[\pi ]$ not dividing $p$, and hence that the same is true for $%
\pi $. For any prime $v$ dividing $p$ of a field containing $\pi $, we let 
\begin{equation*}
s_{\pi }(v)=\frac{\ord_{v}(\pi )}{\ord_{v}(q)};
\end{equation*}
thus $s_{\pi }(v)+s_{\pi }(\iota v)=wt(\pi )$. A Weil $q$-number is
determined up to a root of $1$ (as an element of an algebraic number field)
by the numbers $s_{\pi }(v)$ because they determine all of its valuations.
We call $s_{\pi }$ the \emph{slope function} of $\pi $. A Weil $q$-number
that is itself an integer is called a \emph{Weil }$q$\emph{-integer.}

\subsubsection{Weil germs}

Let $\pi $ be a Weil $p^{n}$-number and $\pi ^{\prime }$ a Weil $%
p^{n^{\prime }}$-number in some field. We say $\pi $ and $\pi ^{\prime }$
are \emph{equivalent }if $\pi ^{n^{\prime }}$ and $\pi ^{\prime n}$ differ
by a root of $1$. A \emph{Weil germ }is an equivalence class of Weil
numbers. The weight and slope function of a Weil germ $\pi $ are the weight
and slope function of any representative of it, and $\mathbb{Q}{}[\pi ]$ is
defined to be the smallest field containing a representative of $\pi $. A
Weil germ is determined by its slope function.

Let $W(p^{\infty })$ denote the set of Weil germs represented by elements of 
$\mathbb{Q}{}^{\text{al}}$. It is an abelian group endowed with an action of 
$\Gamma $. Let $W(p^{\infty })_{m,+}$ denote the subset of $W(p^{\infty })$
consisting of Weil germs of weight $m$ represented by algebraic integers;
thus, 
\begin{equation*}
W(p^{\infty })_{m,+}=\{\pi \in W(p^{\infty })\mid s_{\pi }(v)\geq 0\text{, }%
s_{\pi }(v)+s_{\pi }(\iota v)=m\quad \forall v\}.
\end{equation*}

\subsubsection{Classification of abelian varieties}

Let $A_{0}$ be a simple abelian variety over $\mathbb{F}$, and let $A_{1}$
be a model of $A_{0}$ over $\mathbb{F}{}_{q}\subset \mathbb{F}$ with the
property that $\End(A_{1})=\End(A_{0})$. The Frobenius endomorphism $\pi
_{A_{1}}$ of $A_{1}$ is a Weil $q$-integer of weight $1$ in $C(A_{0})$, and
we let $\pi _{A_{0}}$ denote the germ represented by $\pi _{A_{1}}$ --- it
is independent of the choice of $A_{1}/\mathbb{F}_{q}$. The conjugates of $%
\pi _{A_{0}}$ in $\mathbb{Q}{}^{\text{al}}$ form a $\Gamma $-orbit $\Pi
_{A_{0}}$ in $W(p^{\infty })$, and the map $A_{0}\mapsto \Pi _{A_{0}}$ is a
bijection from the set of isomorphism classes of simple abelian varieties
over $\mathbb{F}{}$ onto the set of $\Gamma $-orbits in $W(p^{\infty
})_{1,+} $.

The various invariants of $A_{0}$ can be read off from $\Pi _{A_{0}}$ as
follows. The images of $\mathbb{Q}{}[\pi _{A_{0}}]$ in $\mathbb{Q}{}^{\text{%
al}}$ are the fixed fields of the stabilizers of the different elements of $%
\Pi _{A_{0}}$, and so $[\mathbb{Q}{}[\pi _{A_{0}}]:\mathbb{Q}{}]=|\Pi |$.
The division algebra $D=_{\text{df}}\End^{0}(A_{0})$ has centre $\mathbb{%
\mathbb{Q}{}}[\pi _{A_{0}}]$, and $D$ splits at no real prime of $\mathbb{Q}%
{}[\pi _{A_{0}}]$, splits at each finite prime not dividing $p$, and has
invariant 
\begin{equation*}
\inv_{v}(D)=s_{\pi }(v)[\mathbb{Q}{}[\pi _{A_{0}}]_{v}:\mathbb{Q}%
{}_{p}]\quad \text{mod }\mathbb{Z}{},
\end{equation*}
at each prime $v$ dividing $p$. By class field theory, the order of $D$ in
the Brauer group of $\mathbb{Q}{}[\pi _{A_{0}}]$ is the smallest positive
integer $e$ such that $e\cdot \inv_{v}(D)\in \mathbb{Z}$ for all $v$, and $%
[D\colon \mathbb{Q}{}[\pi _{A_{0}}]]^{\frac{1}{2}}=e$. Moreover, 
\begin{equation*}
2\dim A_{0}=[D:\mathbb{Q}{}[\pi _{A_{0}}]]^{\frac{1}{2}}\cdot \lbrack 
\mathbb{Q}{}[\pi _{A_{0}}]:\mathbb{Q}{}],
\end{equation*}
and so $A_{0}$ has many endomorphisms. The set of slopes of the
Dieudonn\'{e} module of $A_{0}$ is $\{s_{\pi _{A_{0}}}(v)\mid v|p\}$, and an 
$s$ in this set has multiplicity 
\begin{equation*}
\sum_{v\text{, }s_{\pi _{A_{0}}}(v)=s}\frac{2\dim A_{0}\cdot \lbrack \mathbb{%
Q}{}[\pi _{A_{0}}]_{v}:\mathbb{Q}{}_{p}]}{[\mathbb{Q}{}[\pi _{A_{0}}]:%
\mathbb{Q}{}]}.
\end{equation*}

It remains to classify the Weil germs.

\subsubsection{Classification of Weil germs}

Fix a CM-subfield $K$ of $\mathbb{Q}^{\al}$, finite and Galois over $\mathbb{%
Q}$. For a Weil germ $\pi $ in $\mathbb{Q}{}^{\text{al}}$ and a prime $w$ of 
$\mathbb{Q}{}^{\text{al}}$ dividing $p$, let 
\begin{equation*}
f_{\pi }^{K}(w)=s_{\pi }(w)[K_{w}:\mathbb{Q}{}_{p}].
\end{equation*}
Define $W^{K}(p^{\infty })$ to be the set of Weil germs in $\mathbb{Q}{}^{%
\text{al}}$ represented by an element of $K$ and such that $f_{\pi
}^{K}(w)\in \mathbb{Z}{}$. Since $W(p^{\infty })=\bigcup_{K}W^{K}(p^{\infty
})$, it suffices to describe $W^{K}(p^{\infty })$ for each $K$.

Let $F$ be the maximal real subfield of $K$, and let $X$ and $Y$ be the sets
of primes in $K$ and $F$ respectively dividing $p$. Then there is an exact
sequence 
\begin{equation*}
0\rightarrow W^{K}(p^{\infty })\rightarrow \mathbb{Z}{}^{X}\times \mathbb{Z}%
{}\rightarrow \mathbb{Z}{}^{Y}\rightarrow 0.
\end{equation*}
The first map is $\pi \mapsto (f_{\pi }^{K},wt(\pi ))$ and the second is 
\begin{equation*}
(f,m)\mapsto f|Y-n_{0}\cdot m\cdot \sum_{v\in Y}v
\end{equation*}
where $n_{0}=[K_{w}:\mathbb{Q}{}_{p}]$ for any prime $w$ of $K$ dividing $p$
(it is independent of $w$). Thus $\pi \mapsto f_{\pi }^{K}$ identifies $%
W^{K}(p^{\infty })$ with the set of $f\in \mathbb{Z}{}^{X}$ such that $%
f(w)+f(\iota w)=n_{0}\cdot m$ for some integer $m$ (independent of $w$).

Under $A_{0}\leftrightarrow \Pi _{A_{0}}$, the abelian varieties
corresponding to orbits of $\Gamma $ in $W^{K}(p^{\infty })\cap W(p^{\infty
})_{1,+}$ are those with the property that, for every $\sigma \colon \mathbb{%
Q}{}[\pi _{A_{0}}]\hookrightarrow \mathbb{Q}{}^{\text{al}}$, $\sigma \mathbb{%
Q}{}[\pi _{A_{0}}]\subset K$ and $\End^{0}(A_{0})\otimes _{\mathbb{Q}{}[\pi
_{A_{0}}],\sigma }K$ is a matrix algebra. Thus, there is a one-to-one
correspondence between the isogeny classes of abelian varieties over $%
\mathbb{F}{}$ whose endomorphism algebra is split by $K$ in this sense and
the $\Gal(K/\mathbb{Q}{})$-orbits of $f\in \mathbb{Z}{}^{X}$ such that $%
f(w)+f(\iota w)=n_{0}$ and $f(w)\geq 0$ for all $w.$

\begin{remark*}
Given a possible slope function for a Weil germ $\pi $, the Dieudonn\'{e}
module of the corresponding abelian variety imposes restrictions on the
possible factorizations of $p$ in $\mathbb{Q}{}[\pi ]$. For example, suppose
that $A(\pi )$ has slopes $0$, $\frac{1}{2}$, and $1$, and that the
multiplicity of $\frac{1}{2}$ is $2.$ Then the Dieudonn\'{e} module of $%
A(\pi )$ has a simple isogeny factor of rank $2$, which implies that a prime 
$w$ for which $s_{\pi }(w)=\frac{1}{2}$ must be of degree $2$ (if it had
degree $1$, the action of $\mathbb{Q}{}[\pi ]_{w}$ on the Dieudonn\'{e}
module would split off an isogeny factor of rank $1$). Thus, the
endomorphism algebra of such an abelian variety is commutative.
\end{remark*}

\subsection{A.7. Calculation of the groups over $\mathbb{F}$}

Let $A_{0}$ be an abelian variety over $\mathbb{F}{}$. Then $A_{0}$ is
isogenous to a product $A_{1}^{s_{1}}\times \cdots \times A_{t}^{s_{t}}$
with the $A_{i}$ simple and pairwise nonisogenous, and $G(A_{0})\cong
G(A_{1}\times \cdots \times A_{t})$ for $G=L$, $M$, or $P$; moreover, $%
L(A_{0})_{0}\cong L(A_{1})_{0}\times \cdots \times L(A_{t})_{0}$. Thus, in
the following, we assume that $A_{0}$ is a product of pairwise nonisogenous
simple abelian varieties. Then $E\overset{\text{df}}{=}C(A_{0})$ is either a
CM-algebra or the product of a CM-algebra with $\mathbb{Q}{}$ --- the second
case occurs when one of the isogeny factors of $A_{0}$ is a supersingular
elliptic curve. The Rosati involution is complex conjugation on each
CM-factor of $E$ and the identity map $\mathbb{Q}{}$.

\subsubsection{The Lefschetz group}

The description of $L(A_{0})$ as a subgroup of $(\mathbb{G}{}_{m})_{E/%
\mathbb{Q}{}}$ in terms of characters is the same as in the complex case.

Thus, the group $L(A_{0})$ is the subgroup of $(\mathbb{G}{}_{m})_{E/\mathbb{%
Q}{}}$ whose character group is 
\begin{equation*}
\frac{\mathbb{Z}{}^{\Sigma _{E}}}{\{g\mid g=\iota g\text{ and }\sum g(\sigma
)=0\}}.
\end{equation*}
The weight map $w\colon \mathbb{\mathbb{G}{}}_{m}\rightarrow L(A_{0})$
corresponds to the map 
\begin{equation*}
\lbrack g]\mapsto wt(g)\overset{\text{df}}{=}\sum_{\sigma \in \Sigma
_{E}}g(\sigma )
\end{equation*}
on characters, and the homomorphism $t\colon L(A_{0})\rightarrow \mathbb{%
\mathbb{G}{}}_{m}$ giving the action of $L(A_{0})$ on the Tate object sends $%
1$ to the element of $X^{\ast }(L(A_{0}))$ represented by $-\sigma -\iota
\sigma $, any $\sigma \in \Sigma _{E}$.

It suffices to describe $L(A_{0})_{0}$ in the case that $A_{0}$ is simple.
When $A_{0}$ is a supersingular elliptic curve, $L(A_{0})_{0}=\mu _{2}$;
otherwise $L(A_{0})_{0}$ is the subgroup of $(\mathbb{G}_{m})_{E/\mathbb{Q}}$
whose character group is 
\begin{equation*}
\frac{\mathbb{Z}^{\Sigma _{E}}}{\{g\mid g=\iota g\}}.
\end{equation*}
The map $\mu _{2}\rightarrow L(A)_{0}$ corresponds to the map on characters $%
[g]\mapsto \sum g(\sigma )\mod2.$

When $A_{0}$ is simple, the map $\sigma \mapsto \sigma (\pi _{A_{0}})\colon
\Sigma _{E}\rightarrow \Pi _{A_{0}}$ is bijective and commutes with the
action of $\Gamma $, and so identifies $L(A_{0})$ with the torus whose
character group is 
\begin{equation*}
\frac{\mathbb{Z}^{\Pi _{A_{0}}}}{\{g\mid g=\iota g\text{ and }\sum g(\pi
)=0\}}.
\end{equation*}

\subsubsection{The group $P(A_{0})$}

By definition, $P(A_{0})\subset L(A_{0})$, and a character $[g]$ of $%
L(A_{0}) $ is trivial on $P(A_{0})$ if and only if $g(\pi _{A_{0}})=1$,
where $g(\pi _{A_{0}})$ is the Weil germ $\prod_{\sigma \in \Sigma
_{E}}(\sigma \pi _{A_{0}})^{g(\sigma )}$. A Weil germ is $1$ if and only if
its slopes are all zero, and $[g]$ is trivial on $P(A_{0})$ if and only if 
\begin{equation*}
\sum_{\sigma \in \Sigma _{E}}g(\sigma )s_{\sigma \pi _{A_{0}}}(w)=0\text{,
all }w\text{.}
\end{equation*}
Note that $s_{\sigma \pi _{A_{0}}}(w)=s_{\pi _{A_{0}}}(\sigma ^{-1}w)$.
Similarly, a character $[g]$ of $L(A_{0})_{0}$ is trivial on $L(A_{0})_{0}$
if and only if $g(\pi _{A_{0}}/p^{\frac{1}{2}})=g$ where $p^{\frac{1}{2}}$
also denotes the Weil germ represented by the Weil $p$-number $p^{\frac{1}{2}%
}$.

\subsubsection{The motivic group}

Fix a prime $\ell \in S(A_{0})$ (see Appendix B). Let $\Omega _{\lambda }$
be a finite Galois extension of $\mathbb{Q}{}_{\ell }$ splitting $L(A_{0})$,
and let $\chi \in X^{\ast }(L(A_{0}))$. Then $\chi $ is trivial on $M(A_{0})$
if and only if $H^{2r}(A_{0}^{s},\Omega _{\lambda }(r))_{\chi }$ contains a
nonzero algebraic class for some $r$ and $s$, in which case all the spaces $%
H^{2r}(A_{0}^{s},\Omega _{\lambda }(r))_{\chi }$ consist entirely of
algebraic classes.

\subsubsection{Second description of $P(A_{0})$}

Let $K$ be a CM-subfield of $\mathbb{Q}{}^{\text{al}}$, finite and Galois
over $\mathbb{Q}{}$, and let $P^{K}$ be the torus over $\mathbb{Q}{}$ such
that $X^{\ast }(P^{K})=W^{K}(p^{\infty })$ (as a $\Gamma $-module). Assume
that $K$ is large enough to contain the conjugates of $\mathbb{Q}{}[\pi
_{A_{0}}]$ and to split $\End^{0}(A_{0})$. For any character $\chi $ of $%
L(A_{0})$, $\chi (\pi )\in W^{K}(p^{\infty })$. Thus we have a homomorphism $%
[g]\mapsto \lbrack g(\pi )]\colon X^{\ast }(L(A_{0}))\rightarrow
W^{K}(p^{\infty })$, which clearly commutes with the action of $\Gamma $. It
corresponds to a homomorphism $\rho _{A_{0}}\colon P^{K}(p^{\infty
})\rightarrow L(A_{0})$, whose image is $P(A_{0})$.

\begin{example*}
Let $A_{0}$ be isogenous to a product of elliptic curves, $A_{0}\sim
A_{1}\times \cdots \times A_{t}$, no two of which are isogenous. The centre $%
E$ of the endomorphism algebra of $A_{0}$ is the product $E=\prod E_{i}$ of
the centres of the endomorphism algebras of the $A_{i}$. For each $i$,
choose an embedding $\sigma _{i}\colon E_{i}\hookrightarrow \mathbb{Q}{}^{%
\text{al}}$. A character $g$ of $(\mathbb{G}_{m})_{E/\mathbb{Q}{}}$ is
trivial on $L(A_{0})_{0}$ if and only if, for each $i$ for which $A_{i}$ is
ordinary $g(\sigma _{i})=g(\iota \sigma _{i})$, and for each $i$ (there is
at most one) for which $A_{i}$ is supersingular $2|g(\sigma _{i})$.

Let $\pi \in E$ be a Weil $q$-number representing $\pi _{A_{0}}$, and let $%
\pi =(\pi _{1},\ldots ,\pi _{t})$. Then $g$ is trivial on $P(A_{0})_{0}$ if
and only if $g(\pi ^{N})=q^{N\cdot wt(g)/2}$ for some $N$. The statement
(Spiess 1999, Proposition)

\begin{quote}%
%
Let $\alpha _{1},\ldots ,\alpha _{2m}$ be Weil $q$-numbers of elliptic
curves over $\mathbb{F}_{q}{}$ such that $\alpha _{1}\cdots \alpha
_{2m}=q^{m}$; then, after possibly renumbering the $\alpha _{i}$ and
replacing each $\alpha _{i}$ with $\alpha _{i}^{N}$ for some $N$, $\alpha
_{2j-1}\alpha _{2j}=q$ for $j=1,\ldots ,m.$ 
\end{quote}%
%
implies that this holds only if $g$ is trivial on $L(A_{0})_{0}$. Thus $%
P(A_{0})=L(A_{0})$, and so no product of elliptic curves over $\mathbb{F}{}$
has an exotic Tate class.
\end{example*}

\begin{example*}
Let $A_{0}$ be a simple abelian variety over $\mathbb{F}{}$ and let $\pi $
be its Frobenius germ. Assume that there is a prime $v_{1}$ of degree $1$ of 
$\mathbb{Q}{}[\pi ]$ such that $s_{\pi }(v_{1})=0$, $s_{\pi }(\iota v_{1})=1$%
, and $s_{\pi }(v)=1/2$ for $v\neq v_1$. Let $\pi _{1}$ be a Weil $q$-number
representing $\pi $, and let $g$ be a character of $(\mathbb{G}_{m})_{%
\mathbb{Q}{}[\pi ]/\mathbb{Q}{}}$. For any prime $w$ of $\mathbb{Q}{}^{\text{%
al}}$ dividing $p$, 
\begin{equation*}
\ord_{w}(g(\pi _{1}/q^{\frac{1}{2}}))=\frac{\ord_{v_{1}}q}{2}(-g(\sigma
)+g(\iota \sigma ))
\end{equation*}
where $\sigma $ is the unique embedding of $\mathbb{Q}{}[\pi ]$ such that $%
\sigma ^{-1}w=v_{1}$. Therefore, $g$ is trivial on $P(A_{0})_{0}$ if and
only if $g=\iota g$, i.e., if and only if $g$ is trivial on $L(A_{0})_{0}$.
Thus $P(A_{0})=L(A_{0})$, and no power of $A_{0}$ has an exotic Tate class.
In particular, the Tate conjecture holds for the powers of $A_{0}$.

The abelian varieties of ``K3-type'' of Zarhin are covered by this example
(they are the varieties for which, additionally, $[\mathbb{Q}{}[\pi ]\colon 
\mathbb{Q}{}]=2\dim A_{0}$).
\end{example*}

\begin{example*}
Let $A_{0}$ be a simple abelian variety of dimension $>1$ over $\mathbb{F}{}$
and let $\pi $ be its Frobenius germ. Assume that there is a prime $v_{1}$
of $\mathbb{\mathbb{Q}{}}[\pi ]$ whose decomposition group is $\{1,\iota \}$
for which $s_{\pi }(v_{1})=\frac{1}{2}=$ $s_{\pi }(\iota v_{1})$;
assume moreover that $s_{\pi }(v)=0$ or $1$ for all other primes. Let $\pi
_{1}$ be a Weil $q$-number representing $\pi $, and let $\chi $ be a
character of $X^{\ast }(L(A_{0})_{0})$ that is trivial on $P(A_{0})_{0}$. If 
$\chi =m\chi _{1}$ for some $\chi _{1}\in X^{\ast }(L(A_{0})_{0})$, then $%
\chi _{1}$ is also trivial on $P(A_{0})_{0}$. Thus, we may assume that $\chi 
$ is not divisible in $X^{\ast }(L(A_{0})_{0}).$ Let $g=\sum g(\sigma
)\sigma $ be an element of $\mathbb{Z}{}^{\Sigma _{\mathbb{Q}{}[\pi _{0}]}}$
representing $\chi $ and such that $g(\sigma )\neq 0\Rightarrow g(\iota
\sigma )=0$. For any prime $w$ of $\mathbb{Q}{}^{\text{al}}$ dividing $p$, 
\begin{equation*}
\ord_{w}(g(\pi _{1}))/\ord_{w}(q)\equiv \frac{1}{2}g(\sigma )\mod\mathbb{Z}{}
\end{equation*}

\noindent where $\sigma $ is such that $\sigma ^{-1}w=v_{1}$. Hence $%
g(\sigma )$ is even. As $w$ ranges over the primes dividing $p$, $\sigma $
ranges over the elements of $\Sigma _{\mathbb{Q}{}[\pi _{A_{0}}]}$ for which 
$g(\sigma )\neq 0$. This contradicts the fact that $\chi $ is not divisible.
Hence $\chi =0$, and we see that $P(A_{0})=L(A_{0})$. Hence no power of $%
A_{0}$ has an exotic Tate class.

The ``almost ordinary'' abelian varieties of Lenstra and Zarhin are covered
by this example.
\end{example*}

\subsection{A.8. Reduction of abelian varieties with many endomorphisms: the
fundamental diagram}

Fix a prime $w_{0}$ of $\mathbb{Q}{}^{\text{al}}$ dividing $p$, and let $%
\mathbb{F}{}$ be the residue field. As we noted in \S 1, it follows from the
theory of N\'{e}ron models that an abelian variety $A$ over $\mathbb{Q}{}^{%
\text{al}}$ with many endomorphisms has good reduction at $w_{0}$ to an
abelian variety $A_{0}$ over $\mathbb{F}{}$. We shall explain the map $%
A\mapsto A_{0}$ in terms of the above classifications.

Assume $A$ is isotypic, and let $E$ be a CM-subfield of $\End^{0}(A)$ for
which $H^{1}(A,\mathbb{Q}{})$ is free of rank $1$, and let $\Phi $ be the
CM-type on $E$ defined by its action on $H^{1,0}$. Let $\pi _{A_{0}}$ be the
Weil germ of $A_{0}$ in $E$. We fix an embedding $\rho _{0}\colon
E\hookrightarrow \mathbb{Q}^{\text{al}}$, and explain how to construct $\rho
_{0}(\pi _{A_{0}})$. Let $K$ be a CM-subfield of $\mathbb{Q}{}^{\text{al}}$,
finite and Galois over $\mathbb{Q}{}$, and large enough to contain all
conjugates of $E$. As a subfield of $\mathbb{Q}^{\text{al}}$, $K$ acquires a
prime $w_{0}$. For some $h$, $\frak{P}_{w_{0}}^{h}$ will be principal, say $%
\frak{P}_{w_{0}}^{h}{=(}a)$. Let $\alpha =a^{2n}$ where $n$ is the index of
the unit group of the maximal real subfield of $K$ in the full unit group of 
$K$. Then $\psi _{\rho _{0}}(\alpha )$, where $\psi _{\rho _{0}}$ is the
CM-type on $K$ defined in A.4, is a well-defined Weil $p^{2nhf(\frak{P}%
_{w_{0}}/p)}$-integer of weight $1$ lying in $\rho _{0}E$. Its inverse image
in $E$ represents $\pi _{A_{0}}$.

Assume now that $E$ is a field. The value of the function $s_{\pi _{A_{0}}}$
on a prime $v$ of $E$ dividing $p$ is given by the formula 
\begin{equation*}
s_{\pi _{A_{0}}}(v)=\frac{|\Phi (v)|}{|\Sigma _{E}(v)|}\quad \quad \text{%
(***)}
\end{equation*}
where 
\begin{eqnarray*}
\Sigma _{E}(v) &=&\{\sigma \in \Sigma _{E}\mid v=\sigma ^{-1}w_{0}\} \\
\Phi (v) &=&\Phi \cap \Sigma _{E}(v).
\end{eqnarray*}

Suppose $A$ is simple, and that it corresponds to a $\Gal(K/\mathbb{Q}{})$%
-orbit $\Psi $ in $X^{\ast }(S^{K})$. An element $f\in X^{\ast }(S^{K})$ can
be regarded as a function $f\colon \Sigma _{K}\rightarrow \mathbb{Z}{}$.
Define $\bar{f}$ to be the function $X\rightarrow \mathbb{Z}{}$ such that $%
\bar{f}(w)=\sum_{\tau w_{0}=w}f(\tau ),$ i.e., if $f$ is $\sum f(\tau )\tau $%
, then $\bar{f}$ is $\sum f(\tau )\tau w_{0}$. Then $A_{0}$ is isogenous to
a power a simple abelian variety, which corresponds (as in A.6) to the $\Gal%
(K/\mathbb{Q})$-orbit $\{\bar{f}\mid f\in \Psi \}\subset W^{K}(p^{\infty })$.

Let $K$ be a CM-field, finite and Galois over $\mathbb{Q}{}$, and let $F$ be
the maximal totally real subfield of $K$. If no $p$-adic prime of $F$ splits
in $K$, then $S^K=\mathbb{G}_m$ and the only elements of $W^{K}(p^{\infty })$
are those represented by the Weil $p$-numbers $p^{m/2}$. Otherwise, all the $%
p$-adic primes in $F$ split in $K$, and there is an exact commutative
diagram: 
\begin{equation*}
\begin{CD} 0@>>>X^{*}(S^K)
@>{g}>{\left(\begin{smallmatrix}g\\wt(g)\end{smallmatrix}\right)}>%
\mathbb{Z}^{\Sigma_K}\times\mathbb{Z} @>{\left(\begin{smallmatrix}g\\
m\end{smallmatrix}\right)}>{g|F-m\sum_{\sigma\in\Sigma_F}\sigma}>%
\mathbb{Z}^{\Sigma_F} @>>>0\\
@.@V{g}V{[g(\alpha)]}V@V{\left(\begin{smallmatrix} \tau\\m\end{smallmatrix}
\right)}V{\left(\begin{smallmatrix} \tau{w_0}\\ m\end{smallmatrix}
\right)}V@V{\sigma}V{\sigma{v_0}}V\\
0@>>>W^K(p^{\infty})@>{\pi}>{\left(\begin{smallmatrix} f^{K}_{\pi}\\
wt(\pi)\end{smallmatrix} \right)}>\mathbb{Z}^X\times\mathbb{Z}
@>{\left(\begin{smallmatrix} f\\m\end{smallmatrix}
\right)}>{f|Y-n_{0}{\cdot}m\sum_{v{\in}Y}}v>\mathbb{Z}^Y@>>>0 \end{CD}
\end{equation*}
The element above (or to the left of) an arrow is mapped to the element
below (or to the right) by the arrow. The symbol $g(\alpha )$ denotes $\prod
\sigma (\alpha )^{g(\sigma )}$, $n_{0}=[K_{w_{0}}:\mathbb{Q}{}_{p}]$, and $%
v_{0}$ is the prime on $F$ induced by $w_{0}$.

We saw above that an abelian variety $(A,i)$ of CM-type $(E,\Phi )$ reduces
modulo the prime $w_{0}$ of $\mathbb{Q}{}^{\text{al}}$ to an isotypic
abelian variety $A_{0}$ whose Weil germ is determined by (***). Every simple
abelian variety arises in this way: let $A_{0}$ be a simple abelian variety
over $\mathbb{F}{},$ and let $E$ be a CM-field that can be embedded as a
maximal subfield of $\End^{0}(A)$ containing $\mathbb{Q}{}[\pi _{A_{0}}]$;
algebraic number theory shows that $E$ exists, and it is an elementary
exercise to show that there exist CM-types $\Phi $ on $E$ such that $s_{\pi
_{A_{0}}}$ is given by the formula (***); let $A$ be an abelian variety over 
$\mathbb{Q}{}^{\text{al}}$ of CM-type $(E,\Phi )$; it is uniquely determined
up to isogeny, and $A_{0}$ is isogenous to the reduction of $A$ at $w_{0}$.

Thus, to give a lifting (up to isogeny) of $A_{0}$ to characteristic zero is
to give a CM maximal subfield $E$ of $\End^{0}(A)$ and a CM-type on $E$
satisfying (***).

\section{Numerical Equivalence on Abelian Varieties with Many Endomorphisms}

Let $A$ be an abelian variety of dimension $g$ over an algebraically closed
field. In characteristic zero, two cycles in $\mathcal{Z}^{r}(A)$ are \emph{%
homologically equivalent\/} if their classes in $H^{2r}(A,\mathbb{Q}(r))$
are equal, and in characteristic $p\neq 0$, they are $\ell $\emph{%
-homologically equivalent}, $\ell \neq p$, if their classes in the \'{e}tale
cohomology group $H^{2r}(A,\mathbb{Q}_{\ell }(r))$ are equal. Because of
Poincar\'{e} duality and the compatibility of intersection products with cup
products, homological equivalence implies numerical equivalence. It is
generally conjectured that they coincide.

Part (a) of the following theorem is a special case of a theorem of
Grothendieck (Lieberman 1968, Theorem 4), and part (b) is a theorem of
Clozel (Clozel n.d.). The proof is based on that of Clozel.

\begin{theorem}
\begin{enumerate}
\item  For any abelian variety $A$ with many endomorphisms over an
algebraically closed field $k$ of characteristic zero, homological
equivalence coincides with numerical equivalence on $\mathcal{Z}^r(A)$, all $%
r$.

\item  For any abelian variety $A_0$ over $\mathbb{F}$, there exists a set $%
S $ of primes $\ell$ of density $>0$ (depending on $A_0$) for which $\ell$%
-homological equivalence coincides with numerical equivalence on $\mathcal{Z}%
^r(A_0)$, all $r$.
\end{enumerate}
\end{theorem}

\begin{proof}
In the proof, we ignore Tate twists, i.e., we choose an identification of $%
\mathbb{Q}{}\approx \mathbb{Q}{}(1)$ (or $\mathbb{Q}{}_{\ell }\approx 
\mathbb{Q}{}_{\ell }(1)$).

First consider the characteristic zero case. Choose\footnote{%
For each isotypic isogeny factor $A_{i}$ of $A$, choose a CM-field $E_{i}$
in $\End^{0}(A_{i})$ of degree $2\dim A_{i}$, and let $E=\prod E_{i}$. Write 
$H_{1}(A,\mathbb{Q}{})=E\cdot x_{0}$. For any $c\in E^{\times }$ such that $%
\iota _{E}c=-c$, $ax_{0},bx_{0}\mapsto \func{Tr}_{E/\mathbb{Q}{}}(cab)$ is a
Riemann form on $A$, and we can take $D$ to be any divisor whose class it
is. When $A$, $E$, and $\Phi $ are as in this paragraph, one says that $%
(A,i) $, where $i$ is the inclusion $E\hookrightarrow \End^{0}(A)$, is of 
\emph{CM-type }$(E,\Phi )$.} an \'{e}tale CM-algebra $E\subset \End^{0}(A)$
such that $H^{1}(A,\mathbb{Q})$ is free of rank $1$ as an $E$-module and $E$
is stable under the Rosati involution defined by some ample divisor $D$. The
action of $E\otimes _{\mathbb{Q}}\mathbb{R}$ on $H^{1,0}$ defines a CM-type $%
\Phi $ on $E$. We have $\Hom(E,\mathbb{Q}^{\al})=\Phi \sqcup \bar{\Phi}$.

Let $\Omega $ be the smallest subfield of $\mathbb{Q}^{\al}$ containing $%
\sigma E$ for every homomorphism $\sigma \colon E\rightarrow \mathbb{Q}^{\al%
} $. It is a CM-field, finite and Galois over $\mathbb{Q}$. Let $%
H^{r}(A,\Omega )=H^{r}(A,\mathbb{Q})\otimes \Omega $, and let $%
H^{1}(A)_{\sigma }$ be the subspace of $H^{1}(A,\Omega )$ on which $E$ acts
through $\sigma $. Then $H^{1}(A,\Omega )=\bigoplus_{\sigma \in \Phi \sqcup 
\bar{\Phi}}H^{1}(A)_{\sigma }$ and $H^{1}(A)_{\sigma }$ is one-dimensional.
As $H^{r}(A,\Omega )=\bigwedge_{\Omega }^{r}H^{1}(A,\Omega )$, it follows
that 
\begin{equation*}
H^{r}(A,\Omega )=\bigoplus_{I,J,|I|+|J|=r}H^{r}(A)_{I,J}
\end{equation*}
where $I$ and $J$ are subsets of $\Phi $ and $\iota \Phi $ respectively, and 
$H^{r}(A)_{I,J}=_{\text{df}}H^{r}(A)_{I\sqcup J}$ is the subspace on which $%
e\in E$ acts as $\prod_{\sigma \in I\sqcup J}\sigma e$ --- it is of
dimension $1$ and of Hodge type $(|I|,|J|)$. For $x\in H^{r}(A,\Omega )$,
let $x_{I,J}$ denote the projection of $x$ on $H^{r}(A,\Omega )_{I,J}$.
Because $x\mapsto x_{I,J}$ is multiplication by an idempotent $e_{I,J}$ of $%
E\otimes \Omega $, it sends algebraic classes to algebraic classes.

Let $L$ be the class in $H^{2}(A,\mathbb{Q})$ of the divisor $D$. Because $L$
is algebraic, its isotypic components in $H^{2}(A,\Omega )$ are of type $%
(\sigma ,\iota \sigma )$, $\sigma \in \Sigma _{E}$, and, because $L$ defines
a nondegenerate form on $H_{1}(A,\Omega )$, each such component is nonzero.

For each $\sigma $, choose a nonzero element $\omega _{\sigma }$ of $%
H^{1}(A)_{\sigma }$. Then $(\omega _{\sigma })_{\Phi \sqcup \iota \Phi }$ is
a basis for $H^{1}(A,\Omega )$. We may suppose that the $\omega _{\sigma }$
have been chosen so that the $(\sigma ,\iota \sigma )$ component of $L$ is $%
\omega _{\sigma }\omega _{\iota \sigma }$. Denote $\prod_{\sigma \in
I}\omega _{\sigma }\prod_{\sigma \in J}\omega _{\sigma }$ by $\omega _{I,J}$
--- it is a basis for $H^{r}(A)_{I,J}$. For $i\leq g=\dim A$, 
\begin{equation*}
L^{g-i}=(\sum_{\sigma \in \Sigma _{E}}\omega _{\sigma }\omega _{\iota \sigma
})^{g-i}=\sum_{M}(g-i)!\omega _{M,\iota M}
\end{equation*}
where $M$ runs over the subsets of $\Phi $ with $|M|=g-i$. In particular, $%
\omega _{M,\iota M}$ is algebraic. Moreover, 
\begin{equation*}
L^{g-i}\omega _{IJ}=\sum_{|M|=g-i}(g-i)!\omega _{I\cup M,J\cup \iota M}.
\end{equation*}
Only the subsets $M$ disjoint from both $I$ and $J$ contribute to the sum.

We shall need the following theorem of Lieberman (Kleiman 1968, 2A11, 2.2):

\begin{quote}%
%

Let $\mathcal{A}^{r}$ be the space of algebraic classes in $H^{2r}(A,\mathbb{%
Q}{})$; then for $2r\leq g$, the map $L^{g-2r}\colon \mathcal{A}%
^{r}\rightarrow \mathcal{A}^{g-r}$ is an isomorphism.%
\end{quote}%
%

Suppose $\omega _{IJ}$ is algebraic with $|I|+|J|=2r\leq g$. Let $M=I\cap
\iota J$, so that $I=I_{0}\sqcup M$, $J=J_{0}\sqcup \iota M$, $I_{0}\cap
\iota J_{0}=\emptyset $. We shall prove by induction on $|I\cap \iota J|$
that $\omega _{I_{0},J_{0}}$ is also algebraic. If $|I\cap \iota J|=0$,
there is nothing to prove. If not, $|I\cup \iota J|\leq 2r-1$, and there
exists a subset $M$ of $\Phi $ with $g-2r+1$ elements disjoint from $I\cup
\iota J$. Then $\omega _{I\sqcup M,J\sqcup \iota M}$ is nonzero and
algebraic. By Lieberman's theorem, there exists an $x\in \mathcal{A}^{r-1}$
such that $L^{g-2r+2}x=\omega _{I\sqcup M,J\sqcup \iota M}$. If $\omega
_{I^{\prime },J^{\prime }}$ occurs with nonzero coefficient in $x$, then it
is algebraic. But if $\omega _{I^{\prime },J^{\prime }}$ is chosen so that $%
\omega _{I\sqcup M,J\sqcup \iota M}$ occurs with nonzero coefficient in $%
L^{g-2r+2}\omega _{I^{\prime },J^{\prime }}$, then $I_{0}^{\prime }=I_{0}$, $%
J_{0}^{\prime }=J_{0}$. Since $|I^{\prime }\cap \iota J^{\prime }|=|I\cap
J|-2$, the induction hypothesis shows that $\omega _{I_{0},J_{0}}$ is
algebraic.

We now prove the theorem (in the case of characteristic zero). We have to
show that, for each $r\leq g$, the cup-product pairing 
\begin{equation*}
\mathcal{A}^{r}\times \mathcal{A}^{g-r}\rightarrow \mathbb{Q}{}
\end{equation*}
is nondegenerate. Lieberman's theorem shows that the two spaces have the
same dimension, and so it suffices to show that the left kernel is zero.
Thus, let $x$ be a nonzero element of $\mathcal{A}^{r}$, $r\leq g$, and
suppose $\omega _{I,J}$ occurs with nonzero coefficient in $x$. It suffices
to show that $\omega _{I^{\prime },J^{\prime }}$ is algebraic, where $%
I^{\prime }$ and $J^{\prime }$ are the complements of $I$ and $J$ in $\Phi $
and $\iota \Phi $ respectively. From the last paragraph, we know that $%
\omega _{I,J}=\omega _{I_{0}\sqcup M,J_{0}\sqcup \iota M}$ with $\omega
_{I_{0},J_{0}}$ algebraic and $I_{0}$, $\iota J_{0}$, and $M$ disjoint.
Because $\mathcal{A}^{j}\otimes _{\mathbb{Q}{}}\Omega $ is stable under $\Gal%
(\Omega /\mathbb{Q}{})$, $\iota \omega _{I_{0},J_{0}}=\omega _{\iota
J_{0},\iota I_{0}}$ is algebraic. But $\omega _{I^{\prime },J^{\prime
}}=\omega _{\iota J_{0},\iota I_{0}}\cdot \omega _{N,\iota N}$ where $N$ is
the complement of $I_{0}\sqcup \iota J_{0}\sqcup M$ in $\Phi $, which is
obviously algebraic.

Now consider the case $k=\mathbb{F}$. After possibly replacing $A_{0}$ with
an isogenous variety, we may assume that it lifts to an abelian variety $A$
with many endomorphisms in characteristic zero (see A.8). Let $E$ be a
CM-algebra for $A$ as in the first paragraph of the proof. If $\ell $ is
such that $\iota $ is in the decomposition group of some prime $\lambda $ of 
$\Omega $ dividing $\ell $, then the same argument as in characteristic zero
case applies once one replaces $\mathbb{Q}$ with $\mathbb{Q}_{\ell }$ and $%
\Omega $ with $\Omega _{\lambda }$. The Frobenius density theorem shows that
the set of primes $\ell $ such that $\iota $ is the Frobenius element at a
prime $\lambda $ dividing $\ell $ has density $1/[\Omega :\mathbb{Q}{}]$.
For such a prime $\ell $, $\iota $ is in the decomposition group of $\lambda 
$.
\end{proof}

We strengthen (b) of the theorem by showing that the set $S$ can be chosen
so that $\ell $-homological equivalence coincides with numerical equivalence
on $\mathcal{Z}{}^{r}(A_{0}^{s})$ for all $r$ and $s$.

Let $A_{0}$ be an abelian variety over $\mathbb{F}$, and let $E_{0}$ be the
centre $C(A_{0})$ of $\End^{0}(A_{0})$. Let $\Omega _{0}$ be the composite
of $\mathbb{Q}{}[\sqrt{-p}]$ with the all the fields $\sigma E_{0}$ for $%
\sigma \in \Sigma _{E_{0}}$. Define $S(A_{0})$ to be the set of primes $\ell
\neq p$ such that $\iota $ is contained in the decomposition group of $%
\lambda $ for one (hence every) prime $\lambda $ of $\Omega _{0}$ dividing $%
\ell $. Note that $S(A_{0})$ depends only on the finite set of simple
isogeny factors of $A_{0}$; in particular, $S(A_{0})=S(A_{0}^{s})$.

\begin{proposition}
Statement (b) of the theorem holds with $S=S(A_{0})$.
\end{proposition}

\begin{proof}
Suppose $A_{0}$ is isogenous to $A_{1}^{s_{1}}\times \cdots \times
A_{t}^{s_{t}}$ with the $A_{i}$ simple and nonisogenous in pairs. Assume
initially that none of the $A_{i}$ is a supersingular elliptic curve. Then
each $C(A_{i})$ is a CM-field.

For each $i$, let $D_{i}=\End^{0}(A_{i})$, let $m_{i}=[D_{i}:C(A_{i})]^{%
\frac{1}{2}}$, and let $C(A_{i})_{+}$ be the maximal real subfield of $%
C(A_{i})$. Fix an $\ell \in S(A_{0})$. For each $i$, there exists a field $%
F_{i}$ cyclic of degree $m_{i}$ over $C(A_{i})_{+}$ and such that each real
and $\ell $-adic prime of $C(A_{i})_{+}$ splits in $F_{i}$ and the local
degree at each $p$-adic prime is $m_{i}$ (Artin and Tate 1961, p.\ 105,
Theorem 5). Let $E_{i}=F_{i}\cdot C(A_{i})$. Then $E_{i}$ is a CM-field that
splits $D_{i}$ and can be realized as a subfield of $D_{i}$. Therefore (Tate
1968/69, Th\'{e}or\`{e}me 2), $A_{i}$ is isogenous to the reduction of an
abelian variety $\tilde{A}_{i}$ with $\End^{0}(\tilde{A}_{i})=E_{i}$.

After replacing $A_{0}$ with an isogenous variety, we may suppose that it
lifts to the abelian variety $A=_{\text{df}}\tilde{A}_{1}^{s_{1}}\times
\cdots \times \tilde{A}_{t}^{s_{t}}$. The \'{e}tale algebra $E=_{\text{df}%
}E_{1}^{s_{1}}\times \cdots \times E_{t}^{s_{t}}$ acts on $A$ diagonally,
and satisfies the conditions in the first paragraph of the proof Theorem
B.1. The field $\Omega $ generated by the images of $E$ in $\mathbb{Q}^{\al}$
is $\Omega _{0}\cdot F_{1}\cdots F_{t}$. Because of our choice of the $F_{i}$%
, every $\ell $-adic prime in this field is fixed\footnote{%
Let $F_{0}$ be the maximal totally real subfield of $\Omega _{0}$. The
condition that $\iota $ fixes all $\ell $-adic primes in $\Omega _{0}$ means
that, for each $\ell $-adic prime $v$ of $F_{0}$, $\Omega _{0}\otimes
_{F_{0}}(F_{0})_{v}$ is a field. Because $v$ splits in $F_{1}\cdots F_{t}$,
this property is retained by $\Omega $.} by $\iota $. This completes the
proof of the proposition in this case.

When we add a factor $A_{t+1}^{s_{t+1}}$ to $A_{0}$ with $A_{t+1}$ a
supersingular elliptic curve, $\Omega $ is replaced with $\Omega \cdot
E_{t+1}$ where $E_{t+1}$ can be taken to be any quadratic imaginary field in
which $p$ does not split. If we choose $E_{t+1}$ to be $\mathbb{Q}{}[\sqrt{-p%
}]$, then $\Omega \cdot E_{t+1}=\Omega $, and the preceding argument applies.
\end{proof}

Let $A$ be an abelian variety with many endomorphisms over $\mathbb{Q}{}^{%
\text{al}}$, and let $A_{0}$ be its reduction at the prime $w_{0}$. Fix an $%
\ell \neq p$. Then there are canonical isomorphisms $H^{i}(A,\mathbb{Q}%
_{\ell }(j))\rightarrow H^{i}(A_{0},\mathbb{Q}_{\ell }(j))$ for all $i$ and $%
j$. We say that a cohomology class $\gamma \in H^{2r}(A,\mathbb{Q}{}(r))$ is 
$w_{0}$\emph{-algebraic }if its image $\gamma _{\ell }$ in $H^{2r}(A_{0},%
\mathbb{Q}{}_{\ell }(r))$ is in the $\mathbb{Q}{}$-span of the algebraic
classes on $A_{0}$. Every algebraic class is $w_{0}$-algebraic, but not
every $w_{0}$-algebraic class is algebraic.

\begin{theorem}
For any nonzero $w_{0}$-algebraic class $\alpha $ on $A$, there exists a $%
w_{0}$-algebraic class $\alpha ^{\prime }$ such that $\alpha \cup \alpha
^{\prime }\neq 0$.
\end{theorem}

\begin{proof}
Let $\mathcal{A}^{r}(w_{0})$ be the space of $w_{0}$-algebraic classes in $%
H^{2r}(A,\mathbb{Q}{}(r))$. The proof of the characteristic zero case of the
theorem in A.3 will apply with ``algebraic'' replaced by ``$w_{0}$%
-algebraic'' once we have shown that Lieberman's theorem holds for $\mathcal{%
A}^{r}(w_{0})$: for $2r\leq g$, $L^{g-2r}\colon \mathcal{A}%
^{r}(w_{0})\rightarrow \mathcal{A}^{g-r}(w_{0})$ is an isomorphism.

This map is automatically injective, and so we only have to prove
surjectivity.

Let $\gamma $ be a $w_{0}$-algebraic class in $H^{2g-2r}(A,\mathbb{Q}%
{}(g-r)) $; by assumption, the image $\gamma _{\ell }$ of $\gamma $ in $%
H^{2g-2r}(A_{0},\mathbb{Q}_{\ell }(g-r))$ equals the class $\alpha _{\ell }$
of some $\alpha \in \mathcal{Z}^{g-r}(A_{0})\otimes \mathbb{Q}{}$. There
exists a $\gamma ^{\prime }\in H^{2r}(A,\mathbb{Q}{}(r))$ such that $%
L^{g-2r}\gamma ^{\prime }=\gamma $, and Lieberman's theorem says that there
is an $\alpha ^{\prime }\in \mathcal{Z}^{r}(A_{0})\otimes \mathbb{Q}{}$ such
that the cohomology class of $L^{g-2r}\alpha ^{\prime }$ is $\alpha _{\ell }$%
. The images of $\alpha ^{\prime }$ and $\gamma ^{\prime }$ in $H^{2r}(A_{0},%
\mathbb{Q}{}_{\ell }(r))$ map to $\alpha _{\ell }$ and $\gamma _{\ell }$
respectively under the isomorphism $L^{g-r}\colon H^{2r}(A_{0},\mathbb{Q}%
{}_{\ell }(r))\rightarrow H^{2g-2r}(A_{0},\mathbb{Q}{}_{\ell }(g-r))$. As $%
\alpha _{\ell }=\gamma _{\ell }$, this proves that $\gamma ^{\prime }$ is $%
w_{0}$-algebraic.
\end{proof}

\begin{corollary}
Suppose that the $\ell $-adic cohomology class $c_{\ell }$ of $c\in \mathcal{%
Z}^{r}(A_{0})$ is nonzero. If $c_{\ell }$ is the image of a rational
cohomology class on $A$ (i.e., of an element of $H^{2r}(A,\mathbb{Q}{}(r))$%
), then $c$ is not numerically equivalent to zero.
\end{corollary}

\begin{proof}
Immediate consequence of the theorem and the compatibility of the
cup-product pairings.
\end{proof}

The corollary implies that, if every algebraic class on $A_{0}$ ``lifts'' to
a rational cohomology class in characteristic zero, then $\ell $-adic
homological equivalence on $A_{0}$ coincides with numerical equivalence.

\end{appendix}%
%

\bigskip \centerline{{\bf Bibliography.}}

Artin, E., and Tate, J., Class Field Theory, Harvard, Dept. of Mathematics, 
\textbf{1961}.

Clozel, L., Equivalence num\'erique et \'equivalence cohomologique pour les
vari\'et\'es ab\'eliennes sur les corps finis, preprint, \textbf{n.d.}.

Deligne, P., Hodge cycles on abelian varieties (notes by J.S. Milne). In:
Hodge Cycles, Motives, and Shimura Varieties, Lecture Notes in Math. 900,
Springer, Heidelberg, 9--100, \textbf{1982}.

Deligne, P. and Milne, J.S., Tannakian categories, In: Hodge Cycles,
Motives, and Shimura Varieties, Lecture Notes in Math. 900, Springer,
Heidelberg, 101--228, \textbf{1982}.

van Geemen, B., An introduction to the Hodge conjecture for abelian
varieties. Algebraic cycles and Hodge theory (Torino, 1993), 233--252,
Lecture Notes in Math., 1594, Springer, Berlin, \textbf{1994}.

van Geemen, B., Theta functions and cycles on some abelian fourfolds. Math.
Z. 221, no. 4, 617--631, \textbf{1996}.

Jannsen, U., Motives, numerical equivalence, and semisimplicity, Invent.
Math. 107, 447--459, \textbf{1992}.

Kleiman, S., Algebraic cycles and the Weil conjectures. In: Dix Expos\'{e}s
sur la Cohomologie des Sch\'{e}mas, North-Holland, Amsterdam, 359--386, 
\textbf{1968}.

Lange, H., and Birkenhake, Ch., Complex Abelian Varieties, Springer, \textbf{%
\ 1992}.

Lieberman, David I. Numerical and homological equivalence of algebraic
cycles on Hodge manifolds. Amer. J. Math. 90, 366--74, \textbf{1968}.

Milne, J. S. Values of zeta functions of varieties over finite fields. Amer.
J. Math. 108, 297--360, \textbf{1986}.

Milne, J.S., Lefschetz classes on abelian varieties, Duke Math. J. 96,
639--675, \textbf{1999a}.

Milne, J.S., Lefschetz motives and the Tate conjecture, Compositio Math.
117, 47--81, \textbf{1999b}.

Oort, F., The isogeny class of a CM-type abelian variety is defined over a
finite extension of the prime field. J. Pure Appl. Algebra 3, 399--408, 
\textbf{1973}.

H. Pohlmann, Algebraic cycles on abelian varieties of complex multiplication
type, Ann. of Math. (2) 88, 161--180, \textbf{1968}.

Saavedra Rivano, N., Cat\'{e}gories Tannakiennes, Lecture Notes in Math.
265, Springer, \textbf{1972}.

Schoen, C., Hodge classes on self-products of a variety with an
automorphism, Compositio Math. 65, 3--32, \textbf{1988}.

Serre, Jean-Pierre and Tate, John, Good reduction of abelian varieties, Ann.
of Math. (2) 88, 492--517, \textbf{1968}.

Spiess, Michael, Proof of the Tate conjecture for products of elliptic
curves over finite fields, Math. Ann. 314, 285--290,\textbf{\ 1999}.

Tate, John, Endomorphisms of abelian varieties over finite fields, Invent.
Math. 2, 134--144, \textbf{1966}.

Tate, John, Classes d'isog\'{e}nie des vari\'{e}t\'{e}s ab\'{e}liennes sur
un corps fini (d'apr\`{e}s T. Honda), S\`{e}minaire Bourbaki, 21e ann\'{e}e, 
\textbf{1968/69}, no. 352.

Tate, John, Conjectures on algebraic cycles in $l$-adic cohomology, in
Motives (Seattle, WA, 1991), 71--83, Proc. Sympos. Pure Math., Part 1, Amer.
Math. Soc., Providence, RI, \textbf{1994}.

Weil, A., Abelian varieties and the Hodge ring, \OE uvres Scientifiques,
Vol. III, 421--429, \textbf{1977}.

\end{document}